\documentclass[11pt]{article}

\usepackage[T1]{fontenc}
\usepackage{amsmath, amssymb, mathtools}

\usepackage[numbers]{natbib}

\usepackage{float}
\usepackage{fullpage}
\usepackage{multirow}
\usepackage{hyperref}
\usepackage{bm} 

\usepackage{tikz}
\usetikzlibrary{patterns,decorations.pathreplacing,positioning,calc}
\newcommand{\tikzcircle}[2][red,fill=red]{\tikz[baseline=-0.5ex]\draw[#1,radius=#2] (0,0) circle ;}%

\usepackage{amsthm}
\newtheorem{theorem}{Theorem}
\newtheorem{remark}{Remark}

\title{Multiprecision computations with Schwarz methods}
\author{Michal Outrata\footnote{
Department of Numerical Mathematics, Charles University, Prague, Czech republic, outrata@karlin.mff.cuni.cz.\\
Work supported by the PRIMUS grants PRIMUS/25/SCI/022 and UNCE/24/SCI/005 of Charles University.} \; and Daniel B. Szyld\footnote{Department of Mathematics, Temple University, Philadelphia, PA 19122, USA.}}
\date{}

\begin{document}
\maketitle

\begin{abstract}
We explore and analyze the use of multiprecision arithmetic for several classes of Schwarz
methods and preconditioners, where
the approximate solution of the local problems is performed at
a lower precision, i.e., with fewer digits of accuracy than in 
the underlying (double precision) computation. Conditions for the
appropriate round-off criteria for the lower precision are presented.
It is found experimentally that for the model problems about 5 digits of accuracy
are sufficient to achieve the theoretical restrictions, and thus,
single precision suffices for the local solves. Several numerical experiments illustrate the obtained results.
\end{abstract}

\paragraph*{keywords}
Classical Schwarz methods, multiprecision computation, M-matrices

\paragraph*{AMS}
65F10

\section{Introduction}\label{sec_Intro}

We consider the solution of systems of linear algebraic equations of the form
\begin{equation}\label{eqn_SecIntro_Au_eq_f}
A\mathbf{u} = \mathbf{f}, \quad A\in \mathbb{R}^{N\times N}, \mathbf{f}\in \mathbb{R}^N.
\end{equation}
In particular, we analyze the use of multiprecision arithmetic for three types of
Schwarz methods (see~\cite{gander2008schwarz}): 
(i) the \emph{(damped) additive Schwarz method} ((d)AS, see~\cite{lions1988schwarz,additive1987widlundeltit}), 
(ii) the \emph{restricted additive Schwarz method} (RAS, see~\cite{cai1999restricted}), and 
(iii) the \emph{multiplicative Schwarz method} (MS, see~\cite{lions1988schwarz});
 see~\cite{higham2022mixed} for an overview of multiprecision/mixed-precision algorithms. 
Specifically, we study the solution of the local problems (see precise definitions below) using lower precision.
To the best of our knowledge, this is the first time this approach is analyzed for Schwarz methods.

The underlying idea of Schwarz methods, as a part of the wider family of domain decomposition methods, can be summarized as ``divide and conquer'', where the solution of a large problem is approximated by sub-dividing it into many smaller ones that are computationally less demanding than~\eqref{eqn_SecIntro_Au_eq_f}; these are called the {\em local} problems (or subproblems or subdomain problems);
see Section~\ref{sec_IntroSM} for a detailed discussion. 
For matrices obtained by the discretization of a partial differential equation (PDE), the convergence analysis usually focuses on studying the spectral information of the iteration operator. 
When the method is used as a preconditioner, the convergence analysis usually uses the continuous PDE and its discretization, showing a convergence bound independent of the discretization parameter; see, e.g., \cite{dolean2015introduction,smith2004domain,Toselli:2004:DDM}. 
For the algebraic error analysis of the broader class of stationary iterative methods, we refer the reader to~\cite[Chapter 17]{higham2002accuracy} and the references therein.

For \emph{algebraic Schwarz methods}, where analysis does not take advantage of the 
provenance of the system matrix, we are usually satisfied with  
information about the asymptotic convergence factor of the method (see, e.g.,~\cite{benzi2001algebraic,frommer1999weighted,frommer2001algebraic}), whereas more complete spectral information is often available once we couple the system matrix with further information about its origin, leading to a more complete understanding of the method behavior; 
see, e.g.,~\cite{gander2023fundamental,dolean2009optimized,dolean2015effective, gander2006optimized,ganderiterative,gander2022schwarz}. 
Importantly, these methods are usually used as \emph{preconditioners}, i.e., their convergence is further accelerated using Krylov subspace methods. But in our experience, in order to obtain more insight it is often very useful to \emph{first} study the Schwarz methods (or other domain decomposition methods) \emph{as stand-alone solvers}. Then, based on their analysis, we can obtain an insight into or estimate of the type of performance we can expect when we accelerate these methods with a Krylov subspace method. Moreover, 
in this way we often get additional insights into the weak points of the method, which can be then used to propose an improvement such as a 
coarse space, see, e.g.,~\cite{gander2023fundamental,gander2023piecewise} and also~\cite{sarkis2023convergence}. Schwarz methods as solvers are also fundamental for Schwarz asynchronous iterations; see, e.g., \cite{glusa2020scalable,magoules2017asynchronous}.

The goal of multiprecision algorithms is, in general, to reduce the computation, communication and memory costs by working with some portion of the problem/algorithm in lower precision, e.g., replacing standard double-precision data representation with a single-precision or even half-precision representaions in the computationally most challenging part of the algorithm, see, e.g.,~\cite{higham2022mixed} and the references therein. Intuitively, working (partially) in a lower precision can introduce new issues, e.g., numerical error propagation, and thus introduces a trade-off between computational complexity (by virtue of lowering the precision) and  the level of approximation difficulty (e.g., the floating point precision used). However, as shown in~\cite{carson2018accelerating}, using multiprecision algorithms, where different parts of the algorithm are carried out in different precision, it is sometimes possible to get 
the best of both worlds. To fix ideas, let us consider having a \emph{working} precision (say, double precision) and denote its \emph{unit round-off} by $u_w$ (say, $u_w\approx 10^{-16}$), see~\cite[p.~3]{higham2002accuracy}.
The natural first step is to consider the multiprecision Schwarz methods where the \emph{local} problems are solved in a ``lower precision'', i.e., in a precision with a unit round-off $u_{\ell}$ such that $u_{\ell} > u_w$ (we shall identify the \emph{precision} with its unit round-off, e.g., by ``precision $u_{\ell}$ is lower than $u_w$'' we mean $u_{\ell} > u_w$). Importantly, working in the lower precision $u_{\ell}$ limits not only the precision of the computation but also its \emph{range}, i.e., without careful scaling small/large numbers that we can represent in precision $u_w$ may underflow/overflow in precision $u_{\ell}$, see, e.g.,~\cite{higham2019squeezing}. We comment on the specifics of multiprecision computations relevant to our interest in Section~\ref{sec_IntroMP}.

We note that similar settings have been already considered. In~\cite{giraud2008mixedprecision}, 
where the authors use the (non-overlapping) AS as a preconditioner for the conjugate gradient method (CG),
they use single precision for the preconditioner solve (i.e., running the AS in single precision) and the rest of the CG method in double precision. In~\cite{anzt2019adaptive}, the authors take the (non-overlapping) block-Jacobi as a preconditioner for CG and based on the 1-norm condition number of each of the diagonal blocks they calculate their inverses in half, single, or double precision and then apply these using dense mat-vec products in parallel for each block; similar methodology has been used also in~\cite{tian2025mixed} and tested for practically relevant and challenging 2D and 3D problems. Similarly, in~\cite{schneck2021impact}, the authors study (overlapping, with coarse space) AS as the preconditioner and the effect of using different  precision and data formats for the subdomain matrices (fixed vs.\ floating point precision as well as dense storage of the inverse vs.\ storage of the Cholesky factor) on the performance of preconditioned CG. Note that the focus in all of these papers is on numerical experiments and observations \emph{about the preconditioned CG method}, i.e., the interaction of the domain decomposition method and the different precision choices is present \emph{only implicitly}. 
The analysis focuses on the (often questionable, see~\cite[Section~5.6, Corollary~5.6.7 and onward]{liesen2013krylov}) condition number bound for CG for the preconditioned system and is not interested in the domain decomposition method of choice as a stand-alone solver.

In this paper we focus on multiprecision Schwarz methods with a lower precision $u_{\ell}$ for the local solves, treating both $u_{\ell}$ and the rounding routine as free variables that can and should be chosen so as to preserve or even improve the convergence of the Schwarz method. To that end we propose specific rounding routines, derive sufficient conditions for the convergence of the resulting multiprecision Schwarz methods and numerically demonstrate their effectiveness. Later in the paper, we also consider the effectiveness of these methods as preconditioners, i.e., with the Krylov subspace method acceleration.

Thus, our contribution consists of analyzing Schwarz methods where the local problems are solved with lower precision. Our analysis provides sufficient conditions, and when these conditions are met, one can calculate the minimum number of digits needed in the approximation to the solution of the local problem to obtain overall convergence. Our experiments with multiple types of discretized partial differential equations indicate that our conditions are satisfied with about 5 digits of accuracy. 

The rest of the manuscript is organized as follows: 
in sections~\ref{sec_IntroSM} and~\ref{sec_IntroMP}, we give a brief introduction to (algebraic) Schwarz methods and to multiprecision computations. In section~\ref{sec_AlgebrAnalMPSM}, we introduce and analyze the multiprecision Schwarz methods and demonstrate their performance on several model problems. In Section~\ref{sec_Comp_fpSM_mpSM}, we then explore some additional avenues for the analysis of multiprecision Schwarz methods and we conclude with some remarks in Section~\ref{sec_Conclusion}.


\section{Algebraic Schwarz methods}\label{sec_IntroSM}

Consider a domain $\Omega$ over which we wish to solve a linear differential equation, and assume that
the problem~\eqref{eqn_SecIntro_Au_eq_f} corresponds to a discretization on some grid on the domain $\Omega$; 
this induces the vector space of the unknowns $\overline{W} \cong \mathbb{R}^N$. We can partition these unknowns into 
groups, obtaining vector spaces $\overline{W}_i, i=1,\dotsc ,p$, each of dimension $\overline{N}_i:= \mathrm{dim}(\overline{W}_i) < N$. We require these grouping to be \emph{pairwise disjoint}, i.e., $\overline{W}_i \cap \overline{W}_j = \{0\}$ if $i\neq j$ when embedding these into $\mathbb{R}^N$, and \emph{covering}, i.e.,
\begin{equation*}
\overline{W} = \sum_{i=1}^{p} \overline{W}_i = \left\{ \mathbf{w} \, | \, \mathbf{w}=\sum\limits_{i=1}^p \mathbf{w}_i \; \mathrm{for \; some} \; \mathbf{w}_i \in \overline{W}_i \right\}.
\end{equation*}
This partition is usually done by exploiting a relation which associates the variables of the discrete problem with the underlying geometry of $\Omega$ so that unknowns corresponding to some $\overline{W}_i$ belong to a single, connected subdomain $\overline{\Omega}_i$ of the physical domain $\Omega$; in such case we have $\Omega = \cup_{i=1}^p \overline{\Omega}_i$. We set the restriction operators $\overline{R}_i\, :\, \mathbb{R}^N \rightarrow \mathbb{R}^{\overline{N}_i}$, corresponding to $\overline{N}_i$-by-$N$ zero-one matrices with full row rank $\overline{N}_i$, and obtain the prolongation operators as the transpose of the restrictions, i.e., $\overline{R}_i^T\, :\, \mathbb{R}^{\overline{N}_i} \rightarrow  \mathbb{R}^N$. We assume that the restriction matrices are chosen so that
\begin{equation}\label{eqn_secIntro_SM__Ri_eq_IdZeroPi}
\overline{R}_i = \begin{bmatrix} I_{\overline{N}_i} 0 \end{bmatrix} {\Pi}_{i} \in\mathbb{R}^{\overline{N}_i\times N},
\end{equation}
\noindent where $I_{\overline{N}_i}$ is the identity matrix of the dimension $\overline{N}_i$, $0\in\mathbb{R}^{\overline{N}_i\times (N-\overline{N}_i)}$ and 
${\Pi}_{i}\in\mathbb{R}^{N\times N}$ is a permutation matrix.
Note that
composing the prolongation and restriction we have 
\begin{equation*}
\overline{R}_i^T\overline{R}_i = {\Pi}_i^T \begin{bmatrix} I_{\overline{N}_i} & 0 \\ 0 & 0 \end{bmatrix} {\Pi}_i \in \mathbb{R}^{N\times N}.
\end{equation*}

We also consider a {\em covering} of $\Omega$ with \emph{overlapping} subdomains, by enlarging each subdomain
and obtaining $\Omega_i \supset \bar \Omega_i$, $i=1,\ldots,p$.
The corresponding subspaces are
${W}_i \supset \overline{W}_i$ 
and we set $N_i := \mathrm{dim}(W_i) \geq \overline{N}_i$, with $\sum_{i=1}^p N_i > N$.
The corresponding matrices are 
\begin{equation*}
R_i = \begin{bmatrix} I_{N_i} 0 \end{bmatrix} \Pi_{i}\in\mathbb{R}^{N_i\times N} \in\mathbb{R}^{N_i\times N}
\quad \mathrm{and} \quad
R_i^T R_i = \Pi_i^T \begin{bmatrix} I_{N_i} & 0 \\ 0 & 0 \end{bmatrix} \Pi_i \in \mathbb{R}^{N\times N},
\end{equation*}
\noindent see, e.g.,~\cite{cai1999restricted,efstathiou2003restricted,frommer2001algebraic} for more details. 
Define also the complementary subspaces $\Omega_{\neg i} = \Omega \setminus \Omega_i$, and the corresponding
complementary sets of variables $W_{\neg i}$.

Furthermore, we define the subdomain matrices $A_i$ as the restriction of $A$ to $W_i$, i.e.,
$A_i := R_i A R_i^T$, and similarly, $A_{\neg i} = R_{\neg i} A R_{\neg i}^T$.
We then  have
\begin{equation*}
A = \Pi_i \begin{bmatrix} A_i & K_{i} \\ L_{i} & A_{\neg i} \end{bmatrix} \Pi_i^T \in \mathbb{R}^{N\times N}, \quad \mathrm{for}\; i=1\dotsc ,p.
\end{equation*}
\noindent We then set-up the multi-splitting matrices $M_i$ as
\begin{equation*}
M_i = \Pi_i  \begin{bmatrix} A_i & 0 \\ \times & \ast \end{bmatrix}  \Pi_i^T \in \mathbb{R}^{N\times N}, \quad \mathrm{for}\; i=1\dotsc ,p,
\end{equation*}
\noindent where the blocks $\times$ and $\ast$ can be, in general, chosen arbitrarily but the common choice is to set $\times = 0$ and $\ast = A_{\neg i}$ \cite{frommer1999weighted} or $\ast = \mathrm{diag}(A_{\neg i})$ \cite{benzi2001algebraic,frommer2001algebraic};
for more details on multi-splittings, see the cited works and the references therein.

Equipped with this notation, we can formulate the classical algebraic Schwarz methods -- AS, RAS and MS -- in matrix form
\begin{equation}\label{eqn_SecIntro_SM_IterOperator_T}
\mathbf{u}^{(n+1)} = T_{\star}\mathbf{u}^{(n)} + \mathbf{c}_{\star},
\quad \mathrm{or, \; equivalently, \; for \; the \; errors \;} \quad
\mathbf{e}^{(n+1)} = T_{\star}\mathbf{e}^{(n)},
\end{equation}
\noindent where $\star \in \{ \mathrm{AS}, \mathrm{RAS}, \mathrm{MS} \}$, $\mathbf{c}_{\star} \in \mathbb{R}^{N}$ are some constant vectors, $\mathbf{e}^{(n)} := \mathbf{u} - \mathbf{u}^{(n)} \in \mathbb{R}^{N}$ is the error vector after $n$ iterations and the matrices $T_{\star}$ are the iteration matrices of the respective methods, given by\footnote{An equivalent multi-splitting formulation of these can be found in~\cite[Section 2]{frommer2001algebraic} and~\cite[Sections 2 and 3]{benzi2001algebraic}.}
\begin{equation}\label{eqn_secIntro_MultSchwrz_SchwarzMthds_DefOfIterMtrcsAndPreconds}
\resizebox{.905\textwidth}{!}{$
\begin{gathered}
T_{\mathrm{AS}} := I_N - \sum\limits_{i=1}^{p} R_i^TA_i^{-1} R_i A \equiv I_N - M_{\mathrm{AS}}^{-1}A, \quad 
T_{\mathrm{AS},\theta} := I_N - \theta\sum\limits_{i=1}^{p} R_i^TA_i^{-1} R_i A \equiv I_N - M_{\mathrm{AS},\theta}^{-1}A, \\
T_{\mathrm{RAS}} := I_N - \sum\limits_{i=1}^{p} \overline{R}_i^TA_i^{-1} R_i A \equiv I_N - M_{\mathrm{RAS}}^{-1}A, \quad
T_{\mathrm{MS}} := \prod\limits_{i=p}^{1} \left(I_N -  R_i^TA_i^{-1} R_i A\right) \equiv I_N - M_{\mathrm{MS}}^{-1}A, \quad
\end{gathered}
$}
\end{equation}
\noindent where we also included the \emph{(damped) additive Schwarz method} (dAS) with a damping coefficient $\theta$, 
corresponding to the iteration matrix $T_{\mathrm{AS},\theta}$. Note that we also write each of the iteration matrices 
$T_{\star}$ in the form $I - M_{\star}^{-1}A$, though the notation $M_{\star}^{-1}$ 
does not necessarilly imply that this matrix is nonsinguar.
We use this form of $T_{\star}$
to highlight the fact that the convergence of these stationary methods can be 
further \emph{accelerated} if we reformulate them as \emph{preconditioners} for a Krylov method. 
The preconditioners are then the matrices~$M_{\star}^{-1}$, which for preconditioning we assume to be
nonsingular.
Note that
these preconditioners are to be ``applied'' rather then ``solved with'', i.e., the preconditioner matrix-vector action is given for any vector $\mathbf{v}$ by $\mathbf{v} \mapsto M_{\star}^{-1}\mathbf{v}$. For the additive-based methods, the definition of $M_{\star}^{-1}$ is rather straight-forward while for the multiplicative Schwarz the definition becomes seemingly artificial by having
\begin{equation}\label{eqn_secIntro_MultSchwrz_Tms_into_MmsInv}
M_{\mathrm{MS}}^{-1} =\left( I_N - T_{\mathrm{MS}} \right) A^{-1},
\end{equation} 
\noindent which can be further reformulated for practical use, 
see, e.g.,~\cite[Section 14.3]{saad2003iterative}. 

Based on~\eqref{eqn_SecIntro_SM_IterOperator_T}, we see that convergence (or divergence) for a particular choice of the method is determined by the spectral radius $\rho(T_{\star})$, which has been studied in detail for certain classes of matrices $A$. To that end, 
we denote a  \emph{symmetric, positive-definite} or \emph{SPD} matrix $A$ by  $A \succ 0$.
Recall  that $A$ is a \emph{nonsingular $M$-matrix}, provided that the off-diagonal elements of $A$ are non-positive and all elements of the inverse are non-negative, i.e., $A-\mathrm{diag}(A) \leq 0$ and $A^{-1} \geq 0$, where the inequalities are understood element-wise, see~\cite[Chapter 6]{berman1994nonnegative} or~\cite[Section 2.5]{horn1994topics} for further details and references. 
In the rest of the paper we will simply say {\em an $M$-matrix} meaning a nonsingular $M$-matrix.
We finish this section by recalling convergence results for the classical Schwarz methods for these two classes of matrices,
and where exact arithmetic is assumed.

\begin{theorem}[{\cite[Lemma 2.8]{frommer1999weighted},~\cite[Theorem 3.8]{benzi2001algebraic}}]\label{thm_secIntro_SM_SPDConvergence}
Let $A \succ 0$ and let $q\leq p$ be the smallest number of colors such that we can color all the $p$ 
subspaces $W_1,\dotsc ,W_p$ so that if $W_i\cap W_j \neq \mathbf{0}$, and $i\neq j$ then $W_i$ and $W_j$ have different colors. Then
\begin{equation*}
\rho(T_{MS}) \leq \|T_{MS}\|_A <1
\quad \mathrm{and \; for} \quad
\theta <1/q \quad \mathrm{it \; holds \; that} \quad \rho(T_{AS,\theta}) \leq \|T_{AS,\theta}\|_A <1.
\end{equation*}
\end{theorem}

\begin{theorem}[{\cite[Theorem 3.4]{frommer1999weighted},~\cite[Theorem 4.4]{frommer2001algebraic},~\cite[Theorem 3.5]{benzi2001algebraic}}]\label{thm_secIntro_SM_MmatrixConvergence}
Let $A$ be an $M$-matrix and $q\leq p$ be the smallest number of colors such that we can color all the $p$ subspaces $W_1,\dotsc ,W_p$ so that if $W_i\cap W_j \neq \mathbf{0}$, and $i\neq j$, then $W_i$ and $W_j$ have different colors. Then
\begin{equation*}
\rho(T_{MS}) <1, \quad \rho(T_{RAS}) <1
\quad \mathrm{and \; for} \quad
\theta <1/q \quad \mathrm{it \; holds \; that} \quad \rho(T_{AS,\theta}) <1.
\end{equation*}
\end{theorem}

We remark that there are other methods closely related to the ones mentioned. For example the RAS method has numerous related variants (e.g., WRAS, ASH, RASH, WRASH, see~\cite[Section 6]{frommer2001algebraic} for further references). We do not consider them in this paper.


\section{Multiprecision computations}\label{sec_IntroMP}
As we do \emph{not} work with hardware with a wider selection of precision, the different precisions in our multiprecision algorithms needs to be simulated in some way. The number of options available is limited and both theoretically and practically, two stand out -- the 
\texttt{chop} package
\cite{higham2019simulating} and the 
\texttt{advanpix} package
\cite{advanpix2024}, both implemented in MATLAB. To the best of our knowledge, these are considered the golden standard among the available software for \emph{simulating} various precisions in the numerical analysis and scientific computing community.

\paragraph{\bf advanpix package} Using advanpix, we can specify the number of accurate digits $d_{\ell}$ for each computation, i.e., the package simulates the precision based on the \emph{decadic} notation of numbers in contrast to the binary notation that is commonly used in the hardware, software and also in the IEEE 754-2019 standard for floating point format 
\cite{ieee754_2019}. 
Note that the standard also addresses decimal formats.
Say we want to simulate a \texttt{half} precision (\texttt{fp16}), which corresponds to $u_{\mathrm{half}} \approx 4.88 \times 10^{-4}$. Using \texttt{advanpix}, we have to chose to have either four or five accurate digits, neither of which maps precisely onto the standardized format of \texttt{fp16}. Moreover, \texttt{advanpix} does not include underflow/overflow treatment. 
However, using \texttt{advanpix} allows us to explore also ``new'' precisions which are not yet standardized or even used, e.g., six or eleven accurate digits, and frames the computation precision as more of an ``integer-continuous'' parameter. Moreover, the package is a highly optimized software that overwrites the standard (also highly optimized) MATLAB functions to work with the desired number of accurate digits, e.g., the MATLAB LU or QR factorizations for sparse matrices. Without exploiting these, many problems become too computationally demanding (hence the commercial success of this package). In this context we would also like to highlight that a lot of interest has been recently devoted to \emph{efficient} simulation of arbitrary precisions on GPUs with astonishing results. For example, although the hardware of GPUs is highly optimized only for low-precisions, such as \texttt{fp32}, \texttt{fp16} and even lower, a clever way of simulation of \texttt{fp64} on these GPUs \emph{using these low-precision formats} was competitive with (or even preferable to) the standard hardware implementation of \texttt{fp64}, see~\cite{ozaki2012error,ozaki2013generalization,ozaki2025ozaki} and the references therein. This opens doors to real possibility of efficiently simulating ``new'' low-precisions in practice.

\paragraph*{\bf chop toolbox} The \texttt{chop toolbox} is an open-source MATLAB toolbox\footnote{Towards the end of preparing the manuscript, the \texttt{chop toolbox} has been also released for \texttt{python}, see~\cite{carson2025PyChop}.} developed for simulating different precisions using the native \texttt{double} of MATLAB, essentially by removing a portion of the mantissa of the result after each operation, corresponding to ``rounding'' back to the simulated precision. For computations in \texttt{single} precision or lower, chop faithfully simulates the computation in the precision (see~\cite[Section 3.1]{higham2019simulating}) and can also simulate the underflow/overflow during the computation. Although this toolbox outperforms many other options (see~\cite[Sections 5 and 6]{higham2019simulating}), it makes some computations prohibitively time-consuming, even after adapting it to sparse matrices. Although it allows for arbitrary user-defined formats (defined by the number of bits allocated to the exponent and the significand), we will restrict ourselves to the currently standard ones, summarized in Table~\ref{table_secMultiprecIntro_DiffPrecsSpecs}.

\begin{table}[t]
\begin{center}
\begin{tabular}{c|cccccc}
					& Signif. & Exp. & $u$ & $x_{min}$ & $x_{max}$ \\
\hline
\texttt{q52}		& 5 & 2 & $1.25\times 10^{-1}$ & $6.10\times 10^{-5}$ & $5.73\times 10^{4}$ \\
\texttt{q43}		& 4 & 3 & $6.25\times 10^{-2}$ & $1.56\times 10^{-2}$ & $2.40\times 10^{2}$ \\
\texttt{bfloat16}	& 8 & 8 & $3.91\times 10^{-3}$ & $1.18\times 10^{-38}$ & $3.39\times 10^{38}$ \\
\texttt{fp16}		& 11 & 5 & $4.88\times 10^{-4}$ & $6.10\times 10^{-5}$ & $6.55\times 10^{4}$ \\
\texttt{fp32}		& 24 & 8 & $5.96\times 10^{-8}$ & $1.18\times 10^{-38}$ & $3.40\times 10^{38}$ \\
\texttt{fp64}		& 53 & 11 & $1.11\times 10^{-16}$ & $2.23\times 10^{-308}$ & $1.80\times 10^{308}$ \\
\end{tabular}
\end{center}
\caption{}\label{table_secMultiprecIntro_DiffPrecsSpecs}
\end{table}

In most cases, the computationally most demanding part of Schwarz methods are the subdomain solves, i.e., the operations including $A_i^{-1}$. The issue of underflow/overflow is, in our opinion, an important piece in multiprecision calculations and has been, at least partially, addressed in~\cite{higham2019squeezing}, where the authors propose re-scaling procedures so that (close to) the full range of a given precision is utilized (demonstrated for \texttt{fp16}). To be concrete, having a subdomain problem
\begin{equation}\label{eqn_secMultiprecIntro_Aiui_eq_fi}
A_i \mathbf{u}_i = \mathbf{f}_i,
\end{equation}
(where we omit the iteration index to keep the notation simple) and a precision $u_{\ell}$ with the positive range $[x_{min}^{(\ell)}, x_{max}^{(\ell)}]$, the authors propose several algorithms for calculating and using diagonal matrices for row and column rescaling of~\eqref{eqn_secMultiprecIntro_Aiui_eq_fi}  -- let us denote them $D_i^r$ and $D_i^c$ (corresponding to $R$ and $S$ in~\cite{higham2019squeezing}). For any nonsingular $D_i^r$ and $D_i^c$ we then rewrite~\eqref{eqn_secMultiprecIntro_Aiui_eq_fi} as
\begin{equation}\label{eqn_secMultiprecIntro_Aiui_eq_bi_Rescaled}
\mathcal{A}_i \mathbf{v}_i = \mu \mathbf{b}_i
\end{equation}
\noindent with
\begin{equation*}
\mathcal{A}_i := \mu D_i^r A_i D_i^c,
\quad \mathbf{b}_i := D_i^r \mathbf{f}_i,
\quad \mathbf{u}_i := D_i^c \mathbf{v}_i
\quad \mathrm{and} \quad \mu \in \mathbb{R}.
\end{equation*}
\noindent The goal is to take $D_i^r,D_i^c$ so that $|D_i^r A_i D_i^c| \lesssim 1$ entry-wise and then take $\mu = \nu x_{max}^{(\ell)} $ for some $\nu\in (0,1)$ so that
\begin{equation*}
|\mathcal{A}_i| \equiv |\mu D_i^r A_i D_i^c| \lesssim x_{max}^{(\ell)},
\end{equation*}
\noindent again, in the entry-wise sense. A reasonable choice then is to take $D_i^r$ and $D_i^c$ as in~\cite[Algorithms 2.3 and 2.4]{higham2019squeezing}, i.e.,
to take their diagonal entries 
as the maximum norms of the rows (and then the columns) of $A_i$. According to~\cite[Table 4.5]{higham2019squeezing}, the choice of $\nu = 0.1$ (the authors use $\theta$ in their notation) is reasonable and we comment on this choice later. For the system~\eqref{eqn_secMultiprecIntro_Aiui_eq_bi_Rescaled} we also rescale the right-hand side, namely we write
\begin{equation*}
\mu \mathbf{b}_i = \frac{ \| \mathbf{b}_i \|_{\infty} }{\hat{\nu}_i} \hat{\mathbf{b}}_i
\quad \mathrm{with} \quad
\hat{\mathbf{b}}_i := \hat{\nu}_i \frac{\mu}{\| \mathbf{b}_i \|_{\infty} } \mathbf{b}_i,
\end{equation*}
\noindent where, again, $\hat{\nu}_i \in (0,1)$ allows us to tailor how close to $x_{max}^{(\ell)}$ we rescale the entries of the new right-hand side vector $\hat{\mathbf{b}}_i$. Altogether, we rewrote~\eqref{eqn_secMultiprecIntro_Aiui_eq_fi} into
\begin{equation}\label{eqn_secMultiprecIntro_Aiui_eq_bi_RescaledWithRHS}
\mathcal{A}_i \hat{\mathbf{v}}_i = \hat{\mathbf{b}}_i,
\end{equation}
\noindent which we solve in the precision $u_{\ell}$ and then retrieve $\mathbf{u}_i$ in the precision $u_{w}$ by calculating (also in precision $u_w$)
\begin{equation*}
\mathbf{u}_i = \frac{\hat{\nu}_i}{ \| \mathbf{b}_i \|_{\infty} } D_i^c \hat{\mathbf{v}}_i.
\end{equation*}
\noindent Importantly, the rescaling preserves signs of the entries of the matrix and hence $A_i$ is an $M$-matrix if and only if $\mathcal{A}_i$ is. It can be also adapted to preserve symmetry (see~\cite[Algorithm 2.5]{higham2019squeezing}) and then it also automatically preserves diagonal dominance.

\begin{remark}
The \texttt{advanpix} toolbox is not open source but as per the documentation uses several open source libraries in C and C++. To the best of our understanding, the accuracy of the multiprecision calculations in \texttt{advanpix} uses a large number of bits for the exponent, notably larger than the IEEE 754-2019 standard. Hence, possible low-precision overflow/underflow appearances, e.g., during the LU factorization, are likely prevented by this implementation without the user's knowledge. Hence, the experiments using \texttt{advanpix} should be understood taken with this caveat. 

Nevertheless, we carry out the calculations so as to minimize the chance for underflow, using the appropriate scaling outlined above. For sufficiently small problem sizes we compared the \texttt{advanpix} results with those obtained by \texttt{chop} and saw no significant difference, suggesting no major overflow/underflow occurrences.
\end{remark}


\section{Algebraic analysis of multiprecision Schwarz methods}\label{sec_AlgebrAnalMPSM}
In this section we give analogous results to Theorems~\ref{thm_secIntro_SM_SPDConvergence} and~\ref{thm_secIntro_SM_MmatrixConvergence} when the subdomain solves $A_i^{-1}$ are represented using a lower-precision in some way. The purpose of the numerical experiments here is twofold -- to demonstrate the theoretical results and also to build an intuition for and understanding of the multiprecision Schwarz methods. Therefore, we will use the convergence properties (such as number of iterations or the convergence factor) to compare the results with their ``full precision'' counterparts, as opposed to, e.g., runtimes. All of the code used to produce these is available at \href{https://github.com/MichalOutrata/mpSchwarz}{https://github.com/MichalOutrata/mpSchwarz} but, naturally, the code assumes that both the \texttt{advanpix} as well as the \texttt{chop} toolboxes are available.

We approach the problem from an algebraic point of view, inspired by the results in~\cite{benzi2001algebraic,frommer1999weighted,frommer2001algebraic} assuming exact arithmetic, with the primary goal of carrying out the subdomain solves -- corresponding to $A_i^{-1}$ (or $M_i^{-1}$) -- in a lower precision $u_{\ell}$, compared to the higher working precision $u_w$. Mixed precision is not explicitly mentioned in either of these works, but the sections on inexact solves on the subdomains in ~\cite{benzi2001algebraic,frommer1999weighted,frommer2001algebraic}
can be re-interpreted in terms of lower precision solves.

Following the notation there, we denote with tildes \emph{quantities that have been obtained by precision-reduction in some sense}. For example, if we assume that the matrix $A_i$ is stored in the working precision $u_w$ and we then store it only in a lower precision $u_{\ell}$, the new matrix will be denoted by $\tilde{A}_i$ and 
replacing all ${A}_i$ with $\tilde{A}_i$ in the definition of $\mathcal{A}_i$ or $T_{\star}$ gives us $\tilde{\mathcal{A}}_i$ or 
$\tilde{T}_{\star}$ for $\star \in \{ \mathrm{AS}, \mathrm{RAS}, \mathrm{MS} \}$). 

For clarity, we call the {\em standard rounding procedure}, the rounding to the nearest floating point number as defined in
the IEEE-754 standard \cite{ieee754_2019}.
We emphasize that the symbol $\sim$ does not mean that the quantity was obtained by the standard rounding procedure, quite on the contrary -- we always consider a particular way of obtaining $\tilde{ \mathcal{A} }_i$ from $\mathcal{A}_i$ that suits the situation and is clear from the context. However, we keep a single notation for all of these 
cases (using $\sim$) to highlight the lower-precision nature. We denote the error in the subdomain matrices by $\mathtt{E}_i$, i.e., we have
\begin{equation}\label{eqn_secAlgebraSM_def_mathttEi}
\tilde{A}_i = A_i + \mathtt{E}_i.
\end{equation}
\noindent As is standard in the algebraic convergence theory of Schwarz methods, we are interested in properties of the splittings
\begin{equation}\label{eqn_secAlgebraSM_Ai_eq_tildeAi_min_mathttEi}
A_i = \tilde{A}_i - \left( \tilde{A}_i - A_i \right), \quad i=1,\dotsc ,p \; .
\end{equation}

\subsection{The general case}\label{sec_AlgebraSM_NonSymmCase}
Assuming $A$ is an $M$-matrix, we recall a sufficient condition for the convergence of (damped) AS, RAS and MS is to have
\begin{equation}\label{eqn_secAlgebraSM_tildeAi_inv__tildeAiEi_inv__geq0}
\tilde{A}_i^{-1} \geq 0
\quad \mathrm{and} \quad
\tilde{A}_i^{-1} \left( \tilde{A}_i - A_i \right) = \tilde{A}_i^{-1}\mathtt{E}_i \geq 0.
\end{equation}
\noindent These conditions characterize when the splitting~\eqref{eqn_secAlgebraSM_Ai_eq_tildeAi_min_mathttEi} 
is weak regular (of the first type, see~\cite[Section 4]{frommer2001algebraic}) and thus if~\eqref{eqn_secAlgebraSM_tildeAi_inv__tildeAiEi_inv__geq0} holds for all $i=1,\dotsc ,p$, then (damped) AS, RAS and MS with $A_i^{-1}$ ($M_i^{-1}$) replaced with $\tilde{A}_i^{-1}$ ($\tilde{M}_i^{-1}$) converge, i.e., $\rho(\tilde{T}_{\star}) < 1$. 

In light of the rescaling~\eqref{eqn_secMultiprecIntro_Aiui_eq_fi} to~\eqref{eqn_secMultiprecIntro_Aiui_eq_bi_RescaledWithRHS}, we see that the rounding error is committed at the level of the rescaled system, i.e., instead of solving~\eqref{eqn_secMultiprecIntro_Aiui_eq_bi_RescaledWithRHS} we solve
\begin{equation*}
\tilde{\mathcal{A}}_i \hat{\mathbf{v}}_i = \tilde{\mathbf{b}}_i,
\end{equation*}
\noindent where $\tilde{\mathcal{A}}_i$ (and $\tilde{\mathbf{b}}_i$) is obtained by a rounding technique of our choice applied to $\mathcal{A}_i$ (and $\hat{\mathbf{b}}_i$). In other words, we have
\begin{equation}\label{eqn_secAlgebraSM_tildeAi_eq_muinvDcinvmathcalAiDrinv}
\tilde{A}_i = \mu^{-1}(D_i^r)^{-1} \tilde{\mathcal{A}}_i (D_i^r)^{-1}
\end{equation}
and thus the error matrix $\mathtt{E}_i$ is given by 
\begin{equation*}
\mathtt{E}_i = \tilde{A}_i - A_i = \mu^{-1}(D_i^r)^{-1} \left( \tilde{\mathcal{A}}_i - \mathcal{A}_i \right) (D_i^r)^{-1} = \mu^{-1}(D_i^r)^{-1} \mathtt{F}_i (D_i^r)^{-1},
\end{equation*}
\noindent where we define $\mathtt{F}_i := \tilde{\mathcal{A}}_i - \mathcal{A}_i$ as the rounding error matrix. 

Here we would like to recall a useful observation for diagonal re-scaling of general stationary iterative methods\footnote{We came across this observation in~\cite[Section 17.2, below eq. (17.3)]{higham2002accuracy} but this is likely not the original reference.} -- if the stationary iterative method is based on a splitting $A=M-N$ such that the entries of $M$ are multiples of the corresponding entries of $A$, then the iteration matrices for $A$ and $D_1AD_2$ are similar (i.e., have the same convergence factor) for any nonsingular diagonal matrices $D_1,D_2$. On one hand, this shows that for $u_{\ell}=u_{w}$, the re-scaling above doesn't affect the asymptotic convergence rate in our case. On the other hand, we are clearly interested in cases where the rounding \emph{does} make a difference and through this observation we see that we can expect the convergence factor to be affected by the re-scaling.

Revisiting~\eqref{eqn_secAlgebraSM_tildeAi_inv__tildeAiEi_inv__geq0}, a direct calculation shows that the weak regular splitting conditions are invariant with respect to diagonal scaling with positive entries, i.e., the conditions~\eqref{eqn_secAlgebraSM_tildeAi_inv__tildeAiEi_inv__geq0} are equivalent to
\begin{equation}\label{eqn_secAlgebraSM_WeakRegSplittin_AftrRESCALE}
\tilde{\mathcal{A}}_i^{-1} \geq 0
\quad \mathrm{and} \quad
\tilde{\mathcal{A}}_i^{-1} \left( \tilde{\mathcal{A}}_i - \mathcal{A}_i \right) = \tilde{\mathcal{A}}_i^{-1}\mathtt{F}_i \geq 0.
\end{equation}
\noindent The first ingredient for the analysis of~\eqref{eqn_secAlgebraSM_WeakRegSplittin_AftrRESCALE} is rewritting $\tilde{\mathcal{A}}_i^{-1}$ as
\begin{equation}\label{eqn_secAlgebraSM_tildemathcalAiInvFi_FormulaForNeumannExpansion}
\tilde{\mathcal{A}}_i^{-1} = \mathcal{A}_i^{-1} \left( I + \mathtt{F}_i \mathcal{A}_i^{-1} \right)^{-1},
\end{equation}
\noindent and expanding the inverse matrix there into its Neumann series under the assumption that the spectral radius satisfies
\begin{equation}\label{eqn_secAlgebraSM_mathcalAiInvFi_leq_1}
\rho( \mathcal{A}_i^{-1} \mathtt{F}_i ) < 1,
\end{equation}
Assuming~\eqref{eqn_secAlgebraSM_mathcalAiInvFi_leq_1}, the Neumann serie expansion reads
\begin{equation}\label{eqn_secAlgebraSM_tildemathcalAiInv_NeumannSerPrep}
\tilde{\mathcal{A}}_i^{-1} = \mathcal{A}_i^{-1} \sum\limits_{k=0}^{\infty} \left( -\mathtt{F}_i \mathcal{A}_i^{-1} \right)^k ,
\end{equation}
\noindent and we further rearrange it as
\begin{equation}\label{eqn_secAlgebraSM_tildemathcalAiInv_NeumannSerWithReordr}
\begin{aligned}
\tilde{\mathcal{A}}_i^{-1} &=
  \mathcal{A}_i^{-1} \left( I - \mathtt{F}_i \mathcal{A}_i^{-1} \right) + \mathcal{A}_i^{-1} \left( I - \mathtt{F}_i \mathcal{A}_i^{-1} \right) \mathtt{F}_i \mathcal{A}_i^{-1}\mathtt{F}_i \mathcal{A}_i^{-1} + \dotsc \\
 &= \mathcal{A}_i^{-1} \left( I - \mathtt{F}_i \mathcal{A}_i^{-1} \right) \sum\limits_{k=0}^{\infty} \left( \mathtt{F}_i \mathcal{A}_i^{-1}\right)^{2k} =  \left(\mathcal{A}_i^{-1} - \mathcal{A}_i^{-1}\mathtt{F}_i \mathcal{A}_i^{-1} \right) \sum\limits_{k=0}^{\infty} \left( \mathtt{F}_i \mathcal{A}_i^{-1}\right)^{2k}.
\end{aligned}
\end{equation}

\noindent In order to ensure~\eqref{eqn_secAlgebraSM_WeakRegSplittin_AftrRESCALE} we will focus on ensuring $\mathtt{F}_i\geq 0$ as well as $\tilde{\mathcal{A}}_i^{-1}\geq 0$. Notice that the latter should be natural as we have $A_i^{-1}\geq 0$ and thereby also $\mathcal{A}_i^{-1}\geq 0$, while the condition $\mathtt{F}_i\geq 0$ can be accomplished, at least in theory, by virtue of choosing an appropriate $u_{\ell}$ and the rounding procedure
(see details below). In fact, assuming $\mathtt{F}_i\geq 0$ the natural condition for ensuring also $\tilde{\mathcal{A}}_i^{-1}\geq 0$ (and hence~\eqref{eqn_secAlgebraSM_WeakRegSplittin_AftrRESCALE}) becomes
\begin{equation}\label{eqn_secAlgebraSM_Mmtrx_mathcalAiInv_geq_mathcalAiInvFimathcalAiInv}
\mathcal{A}_i^{-1} \geq \mathcal{A}_i^{-1} \mathtt{F}_i \mathcal{A}_i^{-1},
\end{equation}
\noindent a second condition on the choice of $u_{\ell}$ in addition to~\eqref{eqn_secAlgebraSM_mathcalAiInvFi_leq_1}. Notice that both~\eqref{eqn_secAlgebraSM_mathcalAiInvFi_leq_1} as well as ~\eqref{eqn_secAlgebraSM_Mmtrx_mathcalAiInv_geq_mathcalAiInvFimathcalAiInv} are in some sense generalizations of the standard relative rounding error assumption
\begin{equation}\label{eqn_secAlgebraSM_Mmtrx_mathttFi_leq_usmathcalAi}
| \mathtt{F}_i | \leq u_{\ell} | \mathcal{A}_i |,
\end{equation}
\noindent which is generally guaranteed (the absolute value is to be understood component-wise). Also, similarly to~\cite[Sections 4 and 5]{anzt2019adaptive}, both~\eqref{eqn_secAlgebraSM_mathcalAiInvFi_leq_1} and~\eqref{eqn_secAlgebraSM_Mmtrx_mathcalAiInv_geq_mathcalAiInvFimathcalAiInv} invite us to choose $(u_{\ell})_i$ for each subdomain independently, based on the relevant quantities or their estimates\footnote{In~\cite{anzt2019adaptive}, the authors work with a similar idea but calculate explicitly the analogue of the inverses $\tilde{\mathcal{A}}_i^{-1}$ in different precisions based on their conditioning. This is somewhat complementary to our approach as our interest lies in the analysis of the resulting method rather than in the practical aspect, which has been covered in~\cite{anzt2019adaptive} and we do not comment further on how to choose $u_{\ell}$ (or $(u_{\ell})_i$) for the subdomain problems.}. We keep~\eqref{eqn_secAlgebraSM_mathcalAiInvFi_leq_1} and~\eqref{eqn_secAlgebraSM_Mmtrx_mathcalAiInv_geq_mathcalAiInvFimathcalAiInv} as assumptions, coupling the subdomain problems and the choice of the lower precision $u_{\ell}$ and move our attention to the condition $\mathtt{F}_i \geq 0$. 

If we use the standard round-to-the-nearest rounding as per IEEE 754-2019, then we are unlikely to satisfy $\mathtt{F}_i \geq 0$ except for some special cases. However, the process of rounding is very often fully under our control. Since by definition the off-diagonal entries of $\mathcal{A}_i$ are non-positive while its diagonal entries are non-negative, a simple way to ensure $\mathtt{F}_i \geq 0$ is to take $\tilde{\mathcal{A}}_i$ as the ``sign-informed round up'' of $\mathcal{A}_i$. 

To this end we assume that for any precision $u_{\ell}$ we have at our disposal the functions $\mathtt{rd}_{u_{\ell}}$ and $\mathtt{ru}_{u_{\ell}}$ that round towards zero (down) and towards plus/minus infinity (up)\footnote{In the \texttt{chop} toolbox, these are already implemented and for the  \texttt{advanpix} package, these are straight-forward to implement as we deal with the precision $u_{\ell}$ corresponding to $d_{\ell}$ accurate decimal digits (as opposed to dealing with bits in the case of \texttt{chop}).}. We then introduce the rounding procedure $\mathrm{round}_{Mmtrx}()$ that for any matrix $X$ gives its low-precision approximation $\mathrm{round}_{Mmtrx}(X)$ given by \emph{upward rounding}, i.e.,
\begin{equation*}
\left( \mathrm{round}_{Mmtrx}(X) \right)_{mn} \equiv \left( \mathrm{round}_{Mmtrx}(X,u_{\ell}) \right)_{mn} := 
\begin{dcases}
\mathtt{ru}_{u_{\ell}} \left( \left( X \right)_{mn} \right), \quad \mathrm{if} \; \left( X \right)_{mn} >0, \\
\mathtt{rd}_{u_{\ell}} \left( \left( X \right)_{mn} \right), \quad \mathrm{if} \; \left( X \right)_{mn} <0. \\
\end{dcases}
\end{equation*}
\noindent Taking 
\begin{equation}\label{eqn_secAlgebraSM_Mmtrx_tildemathtcalAi_eq_RoundMmtrx_mathcalAi}
\tilde{\mathcal{A}}_i = \mathrm{round}_{Mmtrx}(\mathcal{A}_i), 
\end{equation}
we get $\mathtt{F}_i \geq 0$ and obtain a convergent multiprecision Schwarz methods under the assumptions \eqref{eqn_secAlgebraSM_mathcalAiInvFi_leq_1} and ~\eqref{eqn_secAlgebraSM_Mmtrx_mathcalAiInv_geq_mathcalAiInvFimathcalAiInv}; we summarize these results in Theorem~\ref{prop_secAlgebraSM_ASRASMSconv_MmtrxCase} below. 

In the following theorem, we use the term {\em full precision} to mean in exact arithemtic, which is when 
Theorems~\ref{thm_secIntro_SM_SPDConvergence} and~\ref{thm_secIntro_SM_MmatrixConvergence} hold.

\begin{theorem}\label{prop_secAlgebraSM_ASRASMSconv_MmtrxCase}
Let $A$ be an $M$-matrix and $q\leq p$ be the smallest number of colors such that we can color all the $p$ subspaces 
$W_1,\dotsc ,W_p$ so that if $W_i\cap W_j \neq \mathbf{0}$, and $i\neq j$, then $W_i$ and $W_j$ have different colors. Moreover, assume that for each $i=1,\dotsc ,p$ 
we replace 
the subdomain solver $A_i^{-1}$ in full precision 
with the subdomain solver $\tilde{A}_i^{-1}$ in a (lower) precision $u_{\ell}$, 
obtaining the multiprecision (damped) AS, RAS and MS methods with the iteration matrices $\tilde{T}_{AS,\theta},\tilde{T}_{RAS}$ 
and $\tilde{T}_{MS}$, respectively. 
Taking $\tilde{A}_i$ as in~\eqref{eqn_secAlgebraSM_tildeAi_eq_muinvDcinvmathcalAiDrinv} with $\tilde{\mathcal{A}}_i$ given as in~\eqref{eqn_secAlgebraSM_Mmtrx_tildemathtcalAi_eq_RoundMmtrx_mathcalAi}, if~\eqref{eqn_secAlgebraSM_mathcalAiInvFi_leq_1} and~\eqref{eqn_secAlgebraSM_Mmtrx_mathcalAiInv_geq_mathcalAiInvFimathcalAiInv} are satisfied, then
\begin{equation*}
\rho(\tilde{T}_{MS}) <1, \quad \rho(\tilde{T}_{RAS}) <1
\quad \mathrm{and \; for} \quad
\theta <1/q \quad \mathrm{it \; holds \; that} \quad  \rho(\tilde{T}_{AS,\theta}) <1,
\end{equation*}
and the multiprecision versions of the classical Schwarz methods are convergent.
\end{theorem}

\begin{remark} \label{Theo:remark}
In practice, in addition to the error $\mathtt{E}_i$ an additional error due to round-off would occur,
say $\mathtt{G}_i$. This very small error may change the sign of the entries. In order to avoid that problem
could explicitely add an additional small fixed perturbation, say $\mathtt{H}_i$, throughout
so that $A_i + \mathtt{E}_i + \mathtt{G}_i + \mathtt{H}_i$ has the desired sign in its entries
(and in the other matrices studied).
In this manner, our analysis would carry over to this practical situation.
We note, however, that in our experiments no such issues were detected.
\end{remark}

Following~\cite[Section 4]{frommer1999weighted},~\cite[Section 7]{frommer2001algebraic} and~\cite[Section 4]{benzi2001algebraic}, we also obtain the comparisons for different choices of $u_{\ell}$. To be more specific, having an $M$-matrix $X$ and two different low-precisions $u_{\ell}^{(1)} \leq u_{\ell}^{(2)}$ with 
$u_{\ell}^{(1)} \leq u_{\ell}^{(2)}$, we obtain 
$
\mathrm{round}_{Mmtrx}(X,u_{\ell}^{(1)}) \leq  \mathrm{round}_{Mmtrx}(X,u_{\ell}^{(2)}),
$
\noindent and hence
\begin{equation*}
\left( \mathrm{round}_{Mmtrx}(X,u_{\ell}^{(1)}) \right)^{-1} \geq  \left( \mathrm{round}_{Mmtrx}(X,u_{\ell}^{(2)})  \right)^{-1}.
\end{equation*}
\noindent Using this for the Schwarz methods, we obtain
\begin{equation}\label{eqn_secAlgebraSM_Mmtrx_ComparisonThmForDiffPrecs}
\rho \left( \tilde{T}_{\star, u_{\ell}^{(1)}} \right) \leq \rho \left( \tilde{T}_{\star, u_{\ell}^{(2)}} \right),
\quad \mathrm{where} \; \star \in \{ \mathrm{AS}, \mathrm{RAS}, \mathrm{MS} \}.
\end{equation}
\noindent In other words, the better precision, the faster convergence. The important questions then become 
\begin{itemize}
\item \textit{Are the conditions~\eqref{eqn_secAlgebraSM_mathcalAiInvFi_leq_1} and~\eqref{eqn_secAlgebraSM_Mmtrx_mathcalAiInv_geq_mathcalAiInvFimathcalAiInv} in some sense sharp or descriptive in the context of convergence of the multiprecision Schwarz methods?}
\item \textit{When is~\eqref{eqn_secAlgebraSM_Mmtrx_ComparisonThmForDiffPrecs} strict?}
\item \textit{Out of those $u_{\ell}$ for which~\eqref{eqn_secAlgebraSM_Mmtrx_ComparisonThmForDiffPrecs} is strict, which should we use, i.e., to what 
extent are there diminishing returns as we approach equality in~\eqref{eqn_secAlgebraSM_Mmtrx_ComparisonThmForDiffPrecs}?}
\item \textit{What computational speedups can we obtain for practically relevant problems using modern HPC architectures?}
\end{itemize}

Next, we investigate the first three questions numerically on model problems coming from discretizations of reaction-advection-diffusion equations. While the last question is itself important, it goes beyond the scope of this manuscript. Regarding the speedups, we note that naive and frequent swapping between rounding modes can be extremely costly as it generally breaks the pipelining of the hardware implementation. However, with modern GPU architectures and with the rise of simulated arithmetic (see the aforementioned Ozaki schemes~\cite{ozaki2025ozaki}), it is plausible that even this issue can be addressed. That being said, 
we do not comment on this further, leaving it open as a possible future direction. In our experiments, the low-precision computations are simulated and hence the slowdown of changing the rounding technique is generally negligible compared to the slowdown due to precision simulations.

Taking the unit square, i.e., $\Omega = \{ \mathbf{x} = [x_1,x_2]^T \in  [0,1]^2 \}$, we consider the partial differential equation
\begin{equation}\label{eqn_secAlgebraSM_Mmtrx_Lu_eq_f}
\mathcal{L}u = f \quad \mathrm{in} \; \Omega 
\quad \mathrm{and} \quad
u = g \quad \mathrm{on} \; \partial \Omega,
\end{equation}
\noindent with the differential operator $\mathcal{L}$ given by
\begin{equation}\label{eqn_secAlgebraSM_Mmtrx_Lu_eq_f_DefinitionOf_L}
\mathcal{L}u := \eta(\mathbf{x})u - \mathrm{div} \left( \alpha(\mathbf{x}) \nabla u \right) + \mathbf{b}(\mathbf{x})\cdot \nabla u.
\end{equation}
\noindent We take the coefficient functions $\eta(\mathbf{x}), \alpha (\mathbf{x})$ and $\mathbf{b}(\mathbf{x}) = [b_1(\mathbf{x}), b_2(\mathbf{x})]^T$ as follows:
\paragraph*{Problem 1} (inspired by~\cite[Figure 2.1]{gander2012optimal})
\begin{equation*}
\eta(\mathbf{x}) := x_1^2 \cos(x_1+x_2)^2, \;\, \alpha(\mathbf{x}) := 20(x_1+x_2)^2 e^{x_1-x_2}
\quad \mathrm{and} \quad
\mathbf{b} (\mathbf{x})  = \left(
\begin{aligned}
 x_2 - 0.5, \\
 x_1 - 0.5
\end{aligned}
\right) \cdot
\end{equation*}
\paragraph*{Problem 2} (inspired by~\cite[Section 4.1]{frommer2023convergence})
\begin{equation*}
\eta \equiv 0, \;\, \alpha \equiv 1 
\quad \mathrm{and} \quad
\mathbf{b} (\mathbf{x})  = \left(
\begin{aligned}
 \beta  (x_1(x_1-1)(1-2x_2)) \\
-\beta (x_2(x_2-1)(1-2x_1))
\end{aligned} 
\right), 
\quad \mathrm{with} \quad \beta = 100.
\end{equation*}
\paragraph*{Problem 3} (based on Problem 2)
\begin{eqnarray*}
\eta \equiv 0, \;\, 
\alpha (\mathbf{x})&=& \left\{
\begin{aligned}
& 10^6 \quad \mathrm{if} \; \|\mathbf{x}-[0.5 \; 0.1]^T\| < 0.25, \\
& 1 \quad \mathrm{otherwise},
\end{aligned} \right. \\
\mathrm{and} \quad
\mathbf{b} (\mathbf{x}) &=& \left(
\begin{aligned}
\beta  (x_1(x_1-1)(1-2x_2)) \\
-\beta (x_2(x_2-1)(1-2x_1))
\end{aligned} 
\right),
\end{eqnarray*}
\noindent again with $\beta = 100$.
To discretize we use the standard 5-point stencil finite difference scheme, adapting some of the code from~\cite{gander2012optimal} and obtain systems of linear equations~\eqref{eqn_SecIntro_Au_eq_f} with $A$ being a non-symmetric $M$-matrix. We then partition $A$ into two overlapping subdomain problems, taking the size of the overlap block to correspond to the bandwidth\footnote{We order the unknowns in the MATLAB's \texttt{numgrid} style, which is column-oriented and results in a banded matrix.} of $A$, i.e., we consider two overlapping subdomains $\Omega_1,\Omega_2 \subset \Omega$ with overlap width\footnote{As a result, we expect the convergence factor of the Schwarz method to deteriorate as $N$ increases, see~\cite[Section 5]{frommer2001algebraic} and \cite[Section 5]{benzi2001algebraic}.} $\mathcal{O}(h)$. We take our right-hand side vector $\mathbf{f}$ and our initial approximation vector $\mathbf{u}^{(0)}$ as random vectors with entries in $(0,1)$.

\begin{figure}[htb]
\centering
\resizebox{1\textwidth}{!}{
	\includegraphics{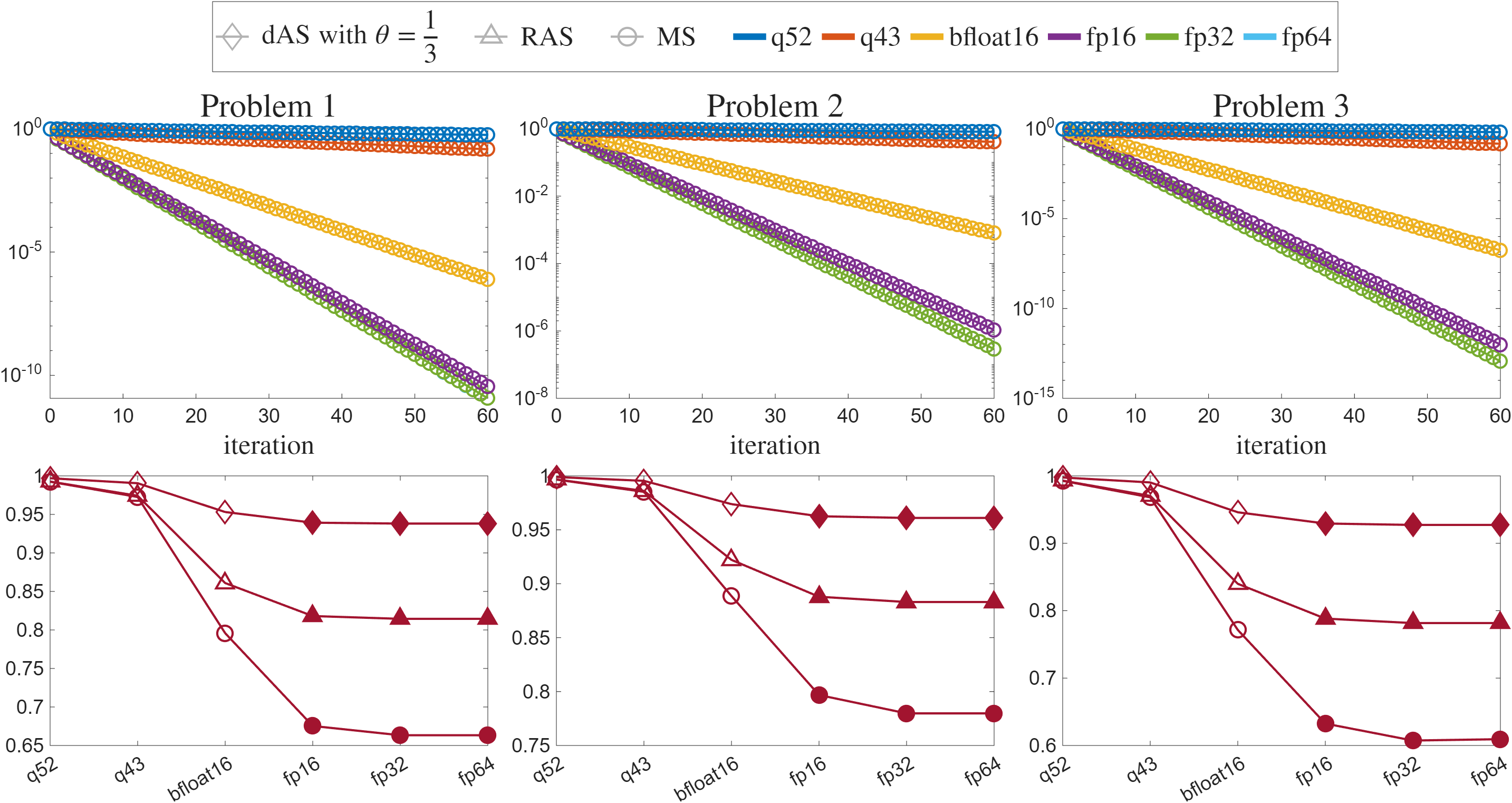}
}
\caption{Top: the 2-norm of the error of the multiplicative Schwarz method with $u_{w}$ corresponding to \texttt{fp64} and different choices of $u_{\ell}$, using the \texttt{chop} toolbox. For \texttt{fp32} and \texttt{fp64} the graphs are indiscernible from each other. Bottom: the observed convergence factor $\rho_{\mathrm{conv}}$ for all Schwarz methods and choices of $u_{\ell}$; if the conditions~\eqref{eqn_secAlgebraSM_mathcalAiInvFi_leq_1} and~\eqref{eqn_secAlgebraSM_Mmtrx_mathcalAiInv_geq_mathcalAiInvFimathcalAiInv} are satisfied for both subdomain $i=1,2$ for a certain $u_{\ell}$, then the marker is filled. For example, for Problem 1 the conditions~\eqref{eqn_secAlgebraSM_mathcalAiInvFi_leq_1} and~\eqref{eqn_secAlgebraSM_Mmtrx_mathcalAiInv_geq_mathcalAiInvFimathcalAiInv} are satisfied for both $i=1,2$ starting from \texttt{fp16}.}\label{fig_secAlgebraSM_NonSymConvCrvs_chop}
\end{figure}

First we fix a $50\times 50$ grid, i.e., $N=2500$, and show the convergence curves and the observed convergence factor\footnote{Having the relative error norms of a Schwarz method, we computed $\rho_{\mathrm{conv}}$ based on a least square fit of linear convergence factor.} $\rho_{\mathrm{conv}}$ in Figure~\ref{fig_secAlgebraSM_NonSymConvCrvs_chop} for the multiplicative Schwarz method and the standard low-precision formats in the \texttt{chop} package (see Table~\ref{table_secMultiprecIntro_DiffPrecsSpecs} above), adjusting the scaling from Section~\ref{sec_IntroMP} so as to use as much of the available range of each precision while not overflowing during the computations. We see that the methods in fact converge in all of the considered precisions, although the conditions~\eqref{eqn_secAlgebraSM_mathcalAiInvFi_leq_1} and~\eqref{eqn_secAlgebraSM_Mmtrx_mathcalAiInv_geq_mathcalAiInvFimathcalAiInv} are satisfied only for \texttt{fp16}, \texttt{fp32} and \texttt{fp64}. Moreover, once the conditions~\eqref{eqn_secAlgebraSM_mathcalAiInvFi_leq_1} and~\eqref{eqn_secAlgebraSM_Mmtrx_mathcalAiInv_geq_mathcalAiInvFimathcalAiInv} are satisfied, they are also satisfied for higher precisions and, more importantly, the observed convergence factor $\rho_{\mathrm{conv}}$ essentially becomes invariant to increasing the precision further. In other words, we get \emph{very little} additional computational benefits (within the first 60 iterations) by considering higher precisions once the conditions~\eqref{eqn_secAlgebraSM_mathcalAiInvFi_leq_1} and~\eqref{eqn_secAlgebraSM_Mmtrx_mathcalAiInv_geq_mathcalAiInvFimathcalAiInv} are satisfied. This suggests that the conditions~\eqref{eqn_secAlgebraSM_mathcalAiInvFi_leq_1} and~\eqref{eqn_secAlgebraSM_Mmtrx_mathcalAiInv_geq_mathcalAiInvFimathcalAiInv} offer a good guidance on the a-priori choice of the working precision $u_w$.

\begin{figure}[t]
\centering
\resizebox{1\textwidth}{!}{
	\includegraphics{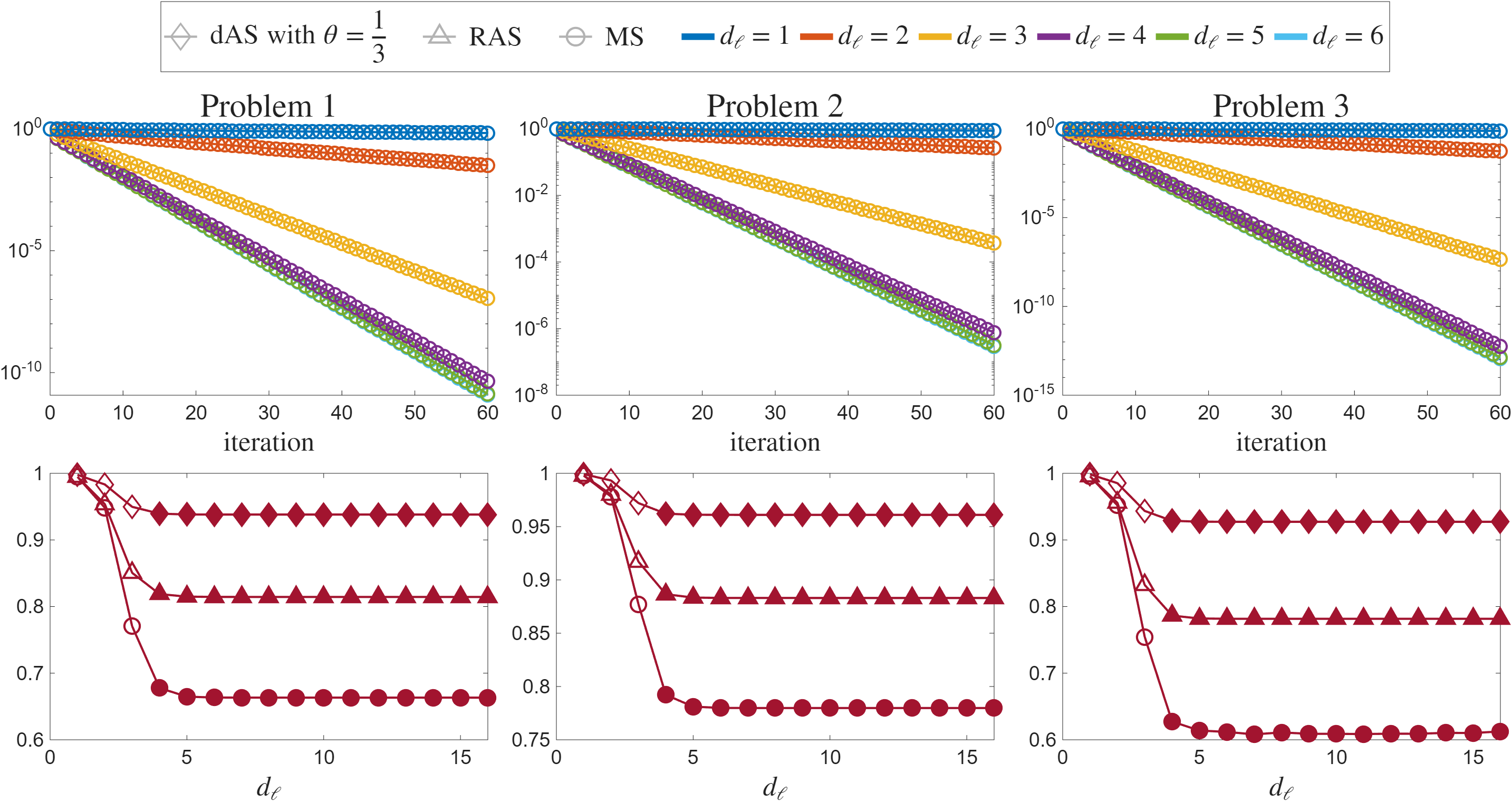}
}
\caption{Top: the 2-norm of the error of the multiplicative Schwarz method with $u_w$ corresponding to \texttt{fp64} and different choices of $d_{\ell}$, using the \texttt{advanpix} toolbox. For $d_{\ell}\geq 6$ the graphs are essentially indiscernible from the green one for $d_{\ell}=5$. Bottom: the observed convergence factor $\rho_{\mathrm{conv}}$ for all Schwarz methods and choices of $d_{\ell}$; if the conditions~\eqref{eqn_secAlgebraSM_mathcalAiInvFi_leq_1} and~\eqref{eqn_secAlgebraSM_Mmtrx_mathcalAiInv_geq_mathcalAiInvFimathcalAiInv} are satisfied for both subdomain $i=1,2$ for a certain $d_{\ell}$, then the marker is filled. For example, for Problem 1 the conditions~\eqref{eqn_secAlgebraSM_mathcalAiInvFi_leq_1} and~\eqref{eqn_secAlgebraSM_Mmtrx_mathcalAiInv_geq_mathcalAiInvFimathcalAiInv} are satisfied for both $i=1,2$ from $d_{\ell} = 4$ onward.}\label{fig_secAlgebraSM_NonSymConvCrvs_advnpx}
\end{figure}

We illustrate this further by showing the analogous experiment but run using the \texttt{advanpix} toolbox, which allows us a finer tuning of the considered precision for the price of giving up the control over underflow/overflow situations (however since we encountered no overflow with \texttt{chop}, this seems not too worrying) in Figure~\ref{fig_secAlgebraSM_NonSymConvCrvs_advnpx}. We see that not only the observations from Figure~\ref{fig_secAlgebraSM_NonSymConvCrvs_chop} still hold true, but the case for the use of the conditions~\eqref{eqn_secAlgebraSM_mathcalAiInvFi_leq_1} and~\eqref{eqn_secAlgebraSM_Mmtrx_mathcalAiInv_geq_mathcalAiInvFimathcalAiInv} as predictors for the suitable precision $u_{\ell}$ is further strengthened.

\begin{remark}\label{rmk_secAlgebraSM_NonSym_NormCondWeakrButSimToEntryCond}
Numerically, the experiments suggest that~\eqref{eqn_secAlgebraSM_mathcalAiInvFi_leq_1} is generally weaker than~\eqref{eqn_secAlgebraSM_Mmtrx_mathcalAiInv_geq_mathcalAiInvFimathcalAiInv}, although we have not been able to establish this as a theoretical result. However, we have never observed this discrepancy to be large, using either \texttt{chop} (e.g., for \texttt{chop} the difference is only present for \texttt{bfloat16} and \texttt{fp16}, where~\eqref{eqn_secAlgebraSM_mathcalAiInvFi_leq_1} was -- for some problems and mesh-sizes -- satisfied for \texttt{bfloat16}, while~\eqref{eqn_secAlgebraSM_Mmtrx_mathcalAiInv_geq_mathcalAiInvFimathcalAiInv} wasn't) or \texttt{advanpix} (again,~\eqref{eqn_secAlgebraSM_mathcalAiInvFi_leq_1} was rarely satisfied for $d_{\ell}=3$, while~\eqref{eqn_secAlgebraSM_Mmtrx_mathcalAiInv_geq_mathcalAiInvFimathcalAiInv} wasn't).
\end{remark}

\begin{figure}[h!]
\centering
\resizebox{0.9\textwidth}{!}{
	\includegraphics{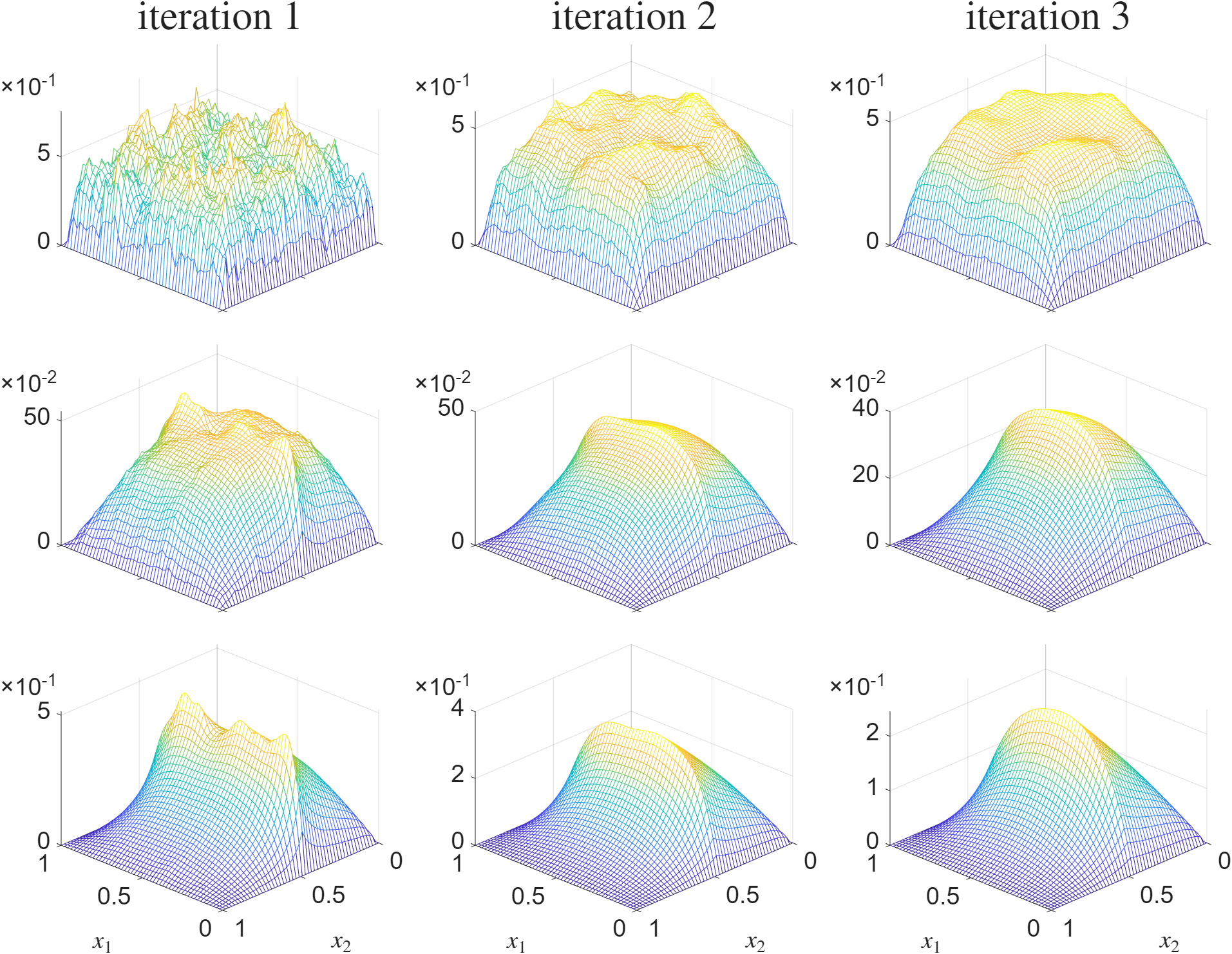}
}
\caption{The errors of the multiplicative Schwarz method used for Problem 1 with $N=2500$ and \texttt{chop} toolbox after 1, 2 and 3 iterations using \texttt{q43} (top), \texttt{bfloat16} and (middle) \texttt{fp16} (bottom). Up to scaling, the error for \texttt{q52} is analogous to the top row and the errors for \texttt{fp32} and \texttt{fp64} are are analogous to the bottom row.}\label{fig_secAlgebraSM_Mmtrx_ErrPlts_chop}
\end{figure}

We further illustrate the tipping point of the conditions~\eqref{eqn_secAlgebraSM_mathcalAiInvFi_leq_1} and~\eqref{eqn_secAlgebraSM_Mmtrx_mathcalAiInv_geq_mathcalAiInvFimathcalAiInv} being met or violated by plotting the error of the multiplicative Schwarz method throughout the initial iterations for different choices of $u_{\ell}$ for Problem 1 in Figure~\ref{fig_secAlgebraSM_Mmtrx_ErrPlts_chop}, again using \texttt{chop} with overflow enabled (but not encountered due to the rescaling). In full precision, we expect to see the classical two-domain profile of the largest eigenmode of the matrix~$T_{\mathrm{MS}}$, which is smooth on each of the subdomains. 
Indeed, for \texttt{fp16}, \texttt{fp32} and \texttt{fp64} that is what we observe. 
However, for \texttt{q52} and \texttt{q43} the ridge in the middle that separates the two subdomains never forms and for \texttt{bfloat16} it takes several iterations to establish to the same extent. In other words, for too low precision $u_{\ell}$ the method effectively looses its continuous level interpretation as a domain decomposition method, although it is still a reasonably effective (even convergent) smoother. Importantly, satisfying the conditions~\eqref{eqn_secAlgebraSM_mathcalAiInvFi_leq_1} and~\eqref{eqn_secAlgebraSM_Mmtrx_mathcalAiInv_geq_mathcalAiInvFimathcalAiInv} is visible not only in the rate of convergence but also in the \emph{nature of it}.

The above observations remained true when changing 
\begin{itemize}
\item the \texttt{chop} and \texttt{advanpix} toolboxes,
\item the problem (we experimented with various settings of reaction-advection-diffusion problems such that $A$ is an $M$-matrix),
\item the method (although, e.g., for dAS we observe an initial period before the error converges to the dominant eigenmode),
\item the initial approximation (the only change is in the initial period before the error converges to the dominant eigenmode).
\end{itemize}

As we kept the problem size relatively small so far, we next experiment also with varying~$N$. However, letting $N$ grow, two numerical limitations come forward -- (i) the \texttt{chop} toolbox becomes too slow and (ii) the verification of the condition~\eqref{eqn_secAlgebraSM_Mmtrx_mathcalAiInv_geq_mathcalAiInvFimathcalAiInv} becomes untenable. Hence, for the following experiments we use only the \texttt{advanpix} toolbox and only verify condition~\eqref{eqn_secAlgebraSM_mathcalAiInvFi_leq_1} (essentially testing whether the observation in Remark~\ref{rmk_secAlgebraSM_NonSym_NormCondWeakrButSimToEntryCond} holds true also for larger $N$).

\begin{figure}[t]
\resizebox{1.\textwidth}{!}{
	\includegraphics{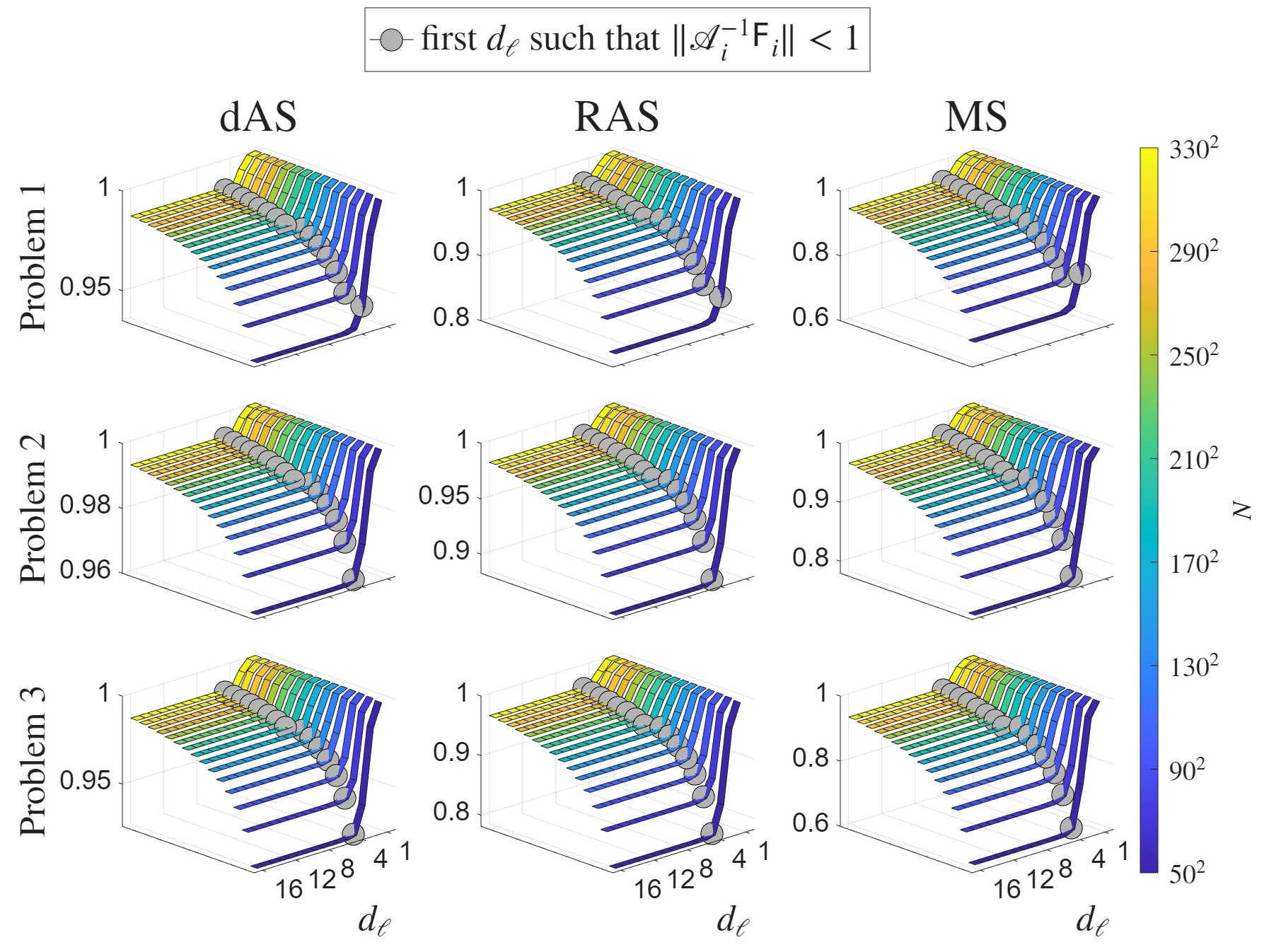}
}
\caption{We show $\rho_{\mathrm{conv}}$ of the Schwarz methods for $d_{\ell}=1,2,\dotsc ,16$ and different problem sizes $N$. The gray circles highlight the first $d_{\ell}$ for which condition~\eqref{eqn_secAlgebraSM_mathcalAiInvFi_leq_1} is satisfied (it is also satisfied for all the following ones).}\label{fig_secAlgebraSM_Mmtrx_ConvFacts}
\end{figure}

All of the above characteristics remained true with the only change being for which $\ell$ the precision $d_{\ell}$ is the first to satisfy the condition~\eqref{eqn_secAlgebraSM_mathcalAiInvFi_leq_1}.
We show these for $N\in \{50^2,\dotsc, \linebreak  330^2\}$ in Figure~\ref{fig_secAlgebraSM_Mmtrx_ConvFacts}. Notably, we see that after the first $d_{\ell}$ such that~\eqref{eqn_secAlgebraSM_mathcalAiInvFi_leq_1} is satisfied there is little to no change in using additional precision, 
precisely as observed before for $N=50^2$. In other words, the dominant eigenmodes of $T_{\star}$ with $\star \in \{ \mathrm{AS}, \mathrm{RAS}, \mathrm{MS} \}$ seem to be well-captured already with limited precision \emph{and} the other eigenmodes are not too sensitive with respect to small perturbations of the subdomain solves and stay ``non-dominant''. In addition, the same type of behavior as shown in Figure~\ref{fig_secAlgebraSM_Mmtrx_ErrPlts_chop} is present for larger $N$, i.e., for too low precision the methods lose their two-domain nature, converge extremely slow but remain effective smoothers.

We see that the convergence factors $\rho_{\mathrm{conv}}$ are remarkably uniform for the different problems as well as with respect to changing the solve precision $u_{\ell}$. Condition~\eqref{eqn_secAlgebraSM_mathcalAiInvFi_leq_1} is satisfied either at $d_{\ell}=4$ or $d_{\ell}=5$, also depending on the size of the problem and once the condition is satisfied, then $\rho_{\mathrm{conv}}$ stabilizes around this final value. In other words, based on these experiments condition~\eqref{eqn_secAlgebraSM_mathcalAiInvFi_leq_1} still governs the required precision. We note that this is perhaps not too surprising as condition~\eqref{eqn_secAlgebraSM_Mmtrx_mathcalAiInv_geq_mathcalAiInvFimathcalAiInv} is clearly \emph{only} a sufficient one -- if it does not hold, then the entries of the matrix $\tilde{\mathcal{A}}_i^{-1}$ in~\eqref{eqn_secAlgebraSM_tildemathcalAiInv_NeumannSerWithReordr} are given as an oscillating sum (rather than a sum of only non-negative numbers), which still can easily sum-up to a non-negative number. On the other hand, if~\eqref{eqn_secAlgebraSM_mathcalAiInvFi_leq_1} doesn't hold, then there's no easy way around it.

We remark that in both~\cite{anzt2019adaptive,schneck2021impact} the authors use the condition numbers of the subdomain matrices for choosing $u_{\ell}$. The condition numbers of the subdomain matrices $A_i$ for both Problem 1 and 2 are fairly small within the range $(10^3,10^5)$, while for Problem 3, the subdomain matrices $A_i$ have condition numbers within the range $(10^9,10^{11})$. We see that the conditioning of the subdomain problem and the precision $u_{\ell}$ seems to interact very little. This remained true for other, similarly focused experiments. We note that the rescaling process does explain part of this observation but we note that the analogous plots for experiments \emph{without the rescaling}, i.e., taking $D_i^{(r,c)} = I_{N_i}$, look fairly similar, although the plots are ``less smooth''.

Notice that as $N$ increases $\rho_{\mathrm{conv}}$ tends towards 1. This is a feature of Schwarz methods -- the convergence factor depends on the physical width of the overlap of the subdomains~$\Omega_1$ and~$\Omega_2$. As noted above, in our setting the overlap width is proportional to $h \sim 1/N$ and therefore this effect is expected even in full precision, which is clearly visible in Figure~\ref{fig_secAlgebraSM_Mmtrx_ConvFacts}.

In addition, note that in condition~\eqref{eqn_secAlgebraSM_mathcalAiInvFi_leq_1} the 2-norm can be replaced by any consistent and equivalent matrix norm and the Neumann series result is still valid; see, e.g.,~\cite[Section 1.3, Lemma 1.3.10]{ortega1990numerical}. In other words, the computationally unfeasible condition~\eqref{eqn_secAlgebraSM_mathcalAiInvFi_leq_1} can be replaced by
\begin{equation*}
\| \mathcal{A}_i^{-1} \mathtt{F}_i \|_F^2 < 1,
\quad \mathrm{or} \quad 
\| \mathcal{A}_i^{-1} \mathtt{F}_i \|_1^2 < 1.
\end{equation*}
\noindent Although these are clearly preferable for the purpose of determining 
the number of digits $d_{\ell}$ (or $d_{\ell_i}, i=1,\dotsc ,p$), see~\cite{anzt2019adaptive}, they might give worse indication of whether or not a given precision is suitable for a given~$N$. 

\begin{figure}[t]
\resizebox{1.\textwidth}{!}{
	\includegraphics{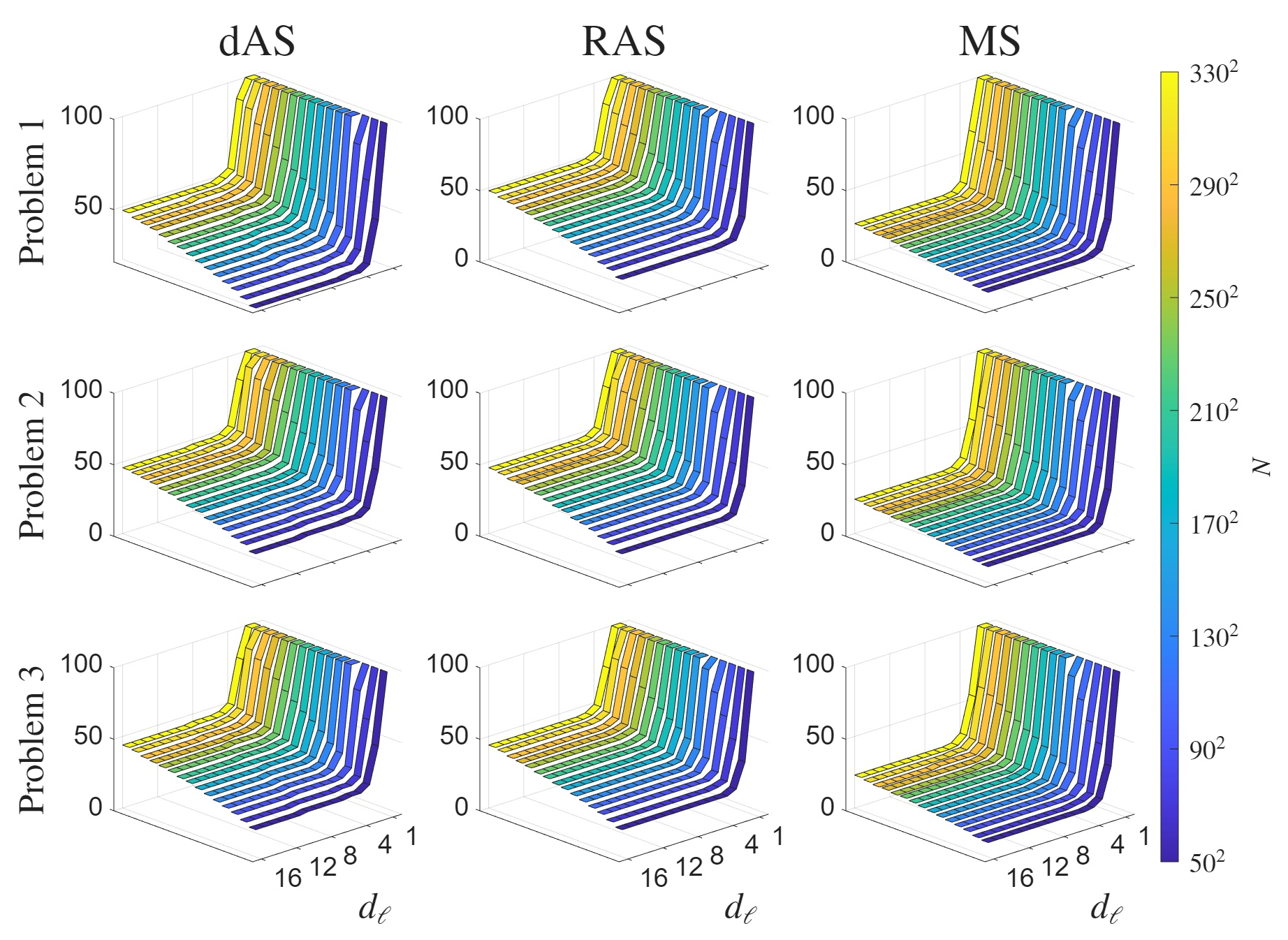}
}
\caption{The number of preconditioned GMRES iterations to reduce the relative residual norm below $10^{-12}$ capped at 100.}\label{fig_secAlgebraSM_Mmtrx_NmbGMRESIter}
\end{figure}

When running Schwarz methods, we can see the effect of the 
lower precision $u_{\ell}$ only on the dominant eigenmode and eigenvalue of $T_{\star}$, as one expects for a fixed-point iteration. However, in practice we usually accelerate Schwarz methods using Krylov subspace methods, i.e., we use Schwarz methods as preconditioners for Krylov methods. In order to be successful preconditioners, calculating in the precision $u_{\ell}$ instead of the precision $u_w$ on the subdomains \emph{should not} make 
the eigenbasis of $M_*^{-1} A$ much more ill-conditioned or the spectrum much more ``spread out'', otherwise a notable slowdown of GMRES convergence (compared to the appropriate full-precision Schwarz method) can occur\footnote{We note that this widely used rule is only a heuristic, see~\cite[Section 5.7]{liesen2013krylov} and the references therein for more details on GMRES convergence analysis.}. In other words, the above experiments do not necessarily suggest that the multiprecision Schwarz methods will be also efficient when used as preconditioners. We investigate that next numerically and use preconditioned GMRES with multiprecision dAS, RAS and MS as the left preconditioners and with preconditioned relative residual tolerance $10^{-12}$, zero initial approximation and maximum number of iterations set to 100. We show the number of GMRES iterations in Figure~\ref{fig_secAlgebraSM_Mmtrx_NmbGMRESIter}.

We observe that the effect of the low-precision does not meaningfully disrupt the number of iterations that preconditioned GMRES needs to reduce the preconditioned relative residual to the tolerance $10^{-12}$. Moreover, we see the same diminishing returns as we did for the convergence factors of the methods in Figure~\ref{fig_secAlgebraSM_Mmtrx_ConvFacts} \emph{and} these occur mostly at the same thresholds, i.e., for the same precisions $u_{\ell}$. We again observe the increase of the iterations as $N$ increases for similar reasons as for the stand-alone solves. While we do not provide the analysis of multiprecision preconditioned GMRES, we refer the reader to~\cite{vieuble2025modular} for an analysis and further references.

\begin{remark}
In the context of (very) low-precision preconditioners, the classical GMRES method should be replaced by the \emph{flexible GMRES} method (FGMRES, see~\cite[Section 9.4.1]{saad2003iterative}) as the low-precision nature of the preconditioner (in our case, the low-precision solves) can be no longer associated with a single matrix representation across the iterations. Moreover, as noted in~\cite{carson2024stability}, FGMRES in the multiprecision context tends to be more robust. In this context, our result show that even for the standard GMRES, not adapted to the specifics of the preconditioners, we obtain the stated results. For more details on FGMRES analysis within the multiprecision context, see~\cite{arioli2008using,carson2024stability}.
\end{remark}

Summarizing, we can say that we observe that in the model examples \emph{four or five digits suffice} to achieve virtually indistinguishable results to full double precision, i.e., running \texttt{fp32} (or even \texttt{fp16} for smaller $N$) should be up to twice as fast (four times as fast) to the standard \texttt{fp64} Schwarz method without any meaningful drawback. Moreover, the results showcase that running \texttt{fp16} or even \texttt{bfloat16} should result in a negligible slowdown while offering up to a further two-fold speed-up.


\subsection{The symmetric case}\label{sec_AlgebraSM_SymmCase}
We consider now \emph{symmetric $M$-matrices}, sometimes also called Stieltjes matrices (Stieltjes matrices are symmetric positive definite, see~\cite[Chapter 6, Theorem 2.3~($D_{16}$)]{berman1994nonnegative}).

First, we would like to highlight that Theorem~\ref{prop_secAlgebraSM_ASRASMSconv_MmtrxCase} applies to this case ``as is''. Moreover, if the rescaling and rounding is done symmetrically (or if we store and work only on, say, the upper-triangular part of the matrix), then the symmetry is preserved on the tilde quantities as well. Notably, to achieve symmetric scaling, the proposed rescaling algorithm needs to be symmetrized, leading to an iterative procedure, see~\cite[Algorithm 2.5]{higham2019squeezing}. Moreover, our rounding routine can be further tailored to preserve other useful properties of the subdomain matrices. 

For example, assuming $A$ has dominant entries on the diagonal, in the sense that $a_{ii} \geq |a_{ij}|$ for all $i,j$,~\cite[Algorithm 2.5]{higham2019squeezing} converges in a single step, yielding $D^{(r)}_i=D^{(c)}_i = \mathrm{diag}(a_{11}^{-1/2}, \dotsc , a_{N_iN_i}^{-1/2})$ with $D^{(r)}_i A_i D^{(c)}_i$ having all ones on the diagonal and the rest of the entries being bounded in absolute value from above by one. Then, taking $\nu$ as some power of two (or other number we represent exactly in precision $u_{\ell}$), we can use the rounding routine
\begin{equation}\label{eqn_secAlgebraSM_tildeAi_ForSPDcase}
\left( \mathrm{round}_{Diag}(X) \right)_{mn} := 
\begin{dcases}
\left( X \right)_{mn} \equiv 1, \quad \mathrm{if} \; m=n, \\
\mathtt{rd}_{u_{\ell}} \left( \left( X \right)_{mn} \right), \quad \mathrm{if} \; m\neq n, \\
\end{dcases}
\end{equation}
\noindent in order to preserve this property (or, e.g., diagonal dominance) also for the rescaled, rounded matrix $\tilde{\mathcal{A}}_i$. Either way, any reasonable rounding should satisfy $\mathtt{F}_i^T = \mathtt{F}_i$ (as $\mathcal{A}_i$ is symmetric), which suffices for now.

Next, we turn our attention to the classical convergence theory for the algebraic Schwarz methods for the symmetric, positive-definite case. As the RAS method is inherently non-symmetric, it is standard to consider the convergence theory only for the (damped) AS and MS methods and as a result we focus only on these two classes\footnote{Some theory for SPD matrices has been developed for variants of the RAS method, see~\cite{cai2003restricted}, and recently, the convergence of RAS for SPD matrices was studied in~\cite{sarkis2023convergence} using variational methods for a simple model problem. In general, convergence of RAS is usually addressed in combination with the particular problem, see~\cite{efstathiou2003restricted}, or based on other properties of the system matrix, see~\cite{frommer2001algebraic}.}. For these, the driving force behind Theorem~\ref{thm_secIntro_SM_SPDConvergence} is the so-called $P$-regular Splitting Theorem, see, e.g.,~\cite[Theorem 7.1.9]{ortega1990numerical}. The core assumption there becomes that the splitting 
in~\eqref{eqn_secAlgebraSM_Ai_eq_tildeAi_min_mathttEi} is a  $P$-regular splitting, i.e., that
\begin{equation}\label{eqn_secAlgebraSM_PregularSplitCond_def}
\tilde{A}_i^T + \tilde{A}_i - A_i \succeq 0.
\end{equation}
\noindent In order to satify~\eqref{eqn_secAlgebraSM_PregularSplitCond_def}, the standard assumption in the literature is $\tilde{A}_i \succeq A_i$ (see, e.g.,~\cite[equation (39), p.621]{benzi2001algebraic}) but, unfortunately, this is not an ``easy to ensure condition'' for a specific rounding routine. Instead, we first observe that for a symmetric scaling, i.e., the case of $D^{(r)}_i=D^{(c)}_i =: D_i$, we obtain 
\begin{equation*}
\tilde{A}_i^T + \tilde{A}_i - A_i  = \mu^{-1} D_i^{-1}\left( \tilde{\mathcal{A}}_i^T + \tilde{\mathcal{A}}_i - \mathcal{A}_i \right) D_i^{-1} = \mu^{-1} D_i^{-1}\left( \tilde{\mathcal{A}}_i + \mathtt{F}_i \right) D_i^{-1},
\end{equation*}
\noindent and ~\eqref{eqn_secAlgebraSM_PregularSplitCond_def} becomes equivalent to 
\begin{equation*}
\tilde{\mathcal{A}}_i + \mathtt{F}_i \succeq 0,
\end{equation*}
\noindent which is ensured by the two following conditions
\begin{equation}\label{eqn_secAlgebraSM_PregularSplitCond_BothOfOurSuffConds}
\tilde{\mathcal{A}}_i \succ 0 
\quad \mathrm{and} \quad
\lambda_{min}( \tilde{\mathcal{A}}_i ) \geq \left| \lambda_{-\infty}( \mathtt{F}_i ) \right| ,
\end{equation}
\noindent where $\lambda_{min}( \tilde{\mathcal{A}}_i ) > 0$ is the smallest eigenvalue of $ \tilde{\mathcal{A}}_i$ (as we assume there $\tilde{\mathcal{A}}_i \succ 0$) and $\lambda_{-\infty}( \mathtt{F}_i )$ is the smallest eigenvalue of $\mathtt{F}_i$ (on the real line, \emph{not} in absolute value, since $\mathtt{F}_i^T = \mathtt{F}_i$). Notice that these conditions differ substantially as we allow for the rounding error matrix to be indefinite.

The first condition in~\eqref{eqn_secAlgebraSM_PregularSplitCond_BothOfOurSuffConds} can be ensured by rounding as in Section~\ref{sec_AlgebraSM_NonSymmCase} so that $\tilde{\mathcal{A}}_i$ is still a Stieltjes matrix and hence symmetric, positive-definite. The second condition can be further expanded on, using the standard perturbation theory of eigenvalues for symmetric matrices (as both $\mathcal{A}_i$ and $\mathtt{F}_i$ are symmetric). Indeed, using Weyl's Theorem (see~\cite[Theorem 4.8 and Corollary 4.9]{stewart1990matrix}) for $ \tilde{\mathcal{A}}_i = \mathcal{A}_i + \mathtt{F}_i$, we obtain
\begin{equation}\label{eqn_secAlgebraSM_PregularSplitCond_WeylThmReslt}
\lambda_{min}( \tilde{\mathcal{A}}_i ) \geq \lambda_{min}( \mathcal{A}_i ) +  \lambda_{-\infty}( \mathtt{F}_i ),
\end{equation}
\noindent so that to ensure the second condition in~\eqref{eqn_secAlgebraSM_PregularSplitCond_BothOfOurSuffConds}, it is enough to require
\begin{equation}\label{eqn_secAlgebraSM_PregularSplitCond_OurSuffConds}
\lambda_{min}( \mathcal{A}_i ) \geq 2 \left| \lambda_{-\infty}( \mathtt{F}_i ) \right|.
\end{equation}
\noindent We summarize the results in the following theorem.

\begin{theorem}\label{prop_secAlgebraSM_ASMSconv_SPDcase}
Let $A$ be a Stieltjes matrix and $q\leq p$ be the smallest number of colors such that we can color all the $p$ subspaces $W_1,\dotsc ,W_p$ so that if $W_i\cap W_j \neq \{0\}$, and $i \neq j$, then $W_i$ and $W_j$ have different colors. Moreover, assume that for each $i=1,\dotsc ,p$ we replace the subdomain solver $A_i^{-1}$ in 
full precision 
with the subdomain solver $\tilde{A}_i^{-1}$ 
with (lower) precision  $u_{\ell}$, 
obtaining the multiprecision (damped) AS and MS methods with the iteration matrices $\tilde{T}_{AS,\theta}$ and $\tilde{T}_{MS}$. Taking $\tilde{A}_i$ as in~\eqref{eqn_secAlgebraSM_tildeAi_eq_muinvDcinvmathcalAiDrinv} with $\tilde{\mathcal{A}}_i$ given as in~\eqref{eqn_secAlgebraSM_Mmtrx_tildemathtcalAi_eq_RoundMmtrx_mathcalAi} with a symmetric scaling, if~\eqref{eqn_secAlgebraSM_mathcalAiInvFi_leq_1},~\eqref{eqn_secAlgebraSM_Mmtrx_mathcalAiInv_geq_mathcalAiInvFimathcalAiInv} and~\eqref{eqn_secAlgebraSM_PregularSplitCond_OurSuffConds} are satisfied, then
\begin{equation*}
\rho(\tilde{T}_{MS}) <1 
\quad \mathrm{and \; for} \quad
\theta <1/q  \quad \mathrm{it \; holds \; that} \quad \rho(\tilde{T}_{AS,\theta}) <1.
\end{equation*}
\noindent and the multiprecision versions of the classical Schwarz methods are convergent.
\end{theorem}

We note that an analysis for SPD matrices and multiprecision additive Schwarz methods has been considered elsewhere; 
see~\cite{anzt2019adaptive,giraud2008mixedprecision,schneck2021impact,tian2025mixed}. However, in all of these papers the authors consider the (non-damped) AS \emph{as a preconditioner for CG} and hence the analysis and/or numerical investigation focus on the preconditioned CG method, e.g., using variational techniques that allow establishing a bound on the condition number of the preconditioned system. We also note that the assumptions look somewhat similar\footnote{Compare~\cite[equations (12) and (14)]{schneck2021impact} and~\cite[Section 5]{anzt2019adaptive} with~\eqref{eqn_secAlgebraSM_mathcalAiInvFi_leq_1} and~\eqref{eqn_secAlgebraSM_PregularSplitCond_OurSuffConds}.}. Indeed, to get a less restrictive version of~\eqref{eqn_secAlgebraSM_PregularSplitCond_OurSuffConds} we can replace $\lambda_{-\infty}( \mathtt{F}_i )$ with $- \max\limits_{\lambda \; \in \; \sigma(\mathtt{F}_i)} |\lambda| \equiv -\rho(\mathtt{F}_i)$ so that instead of~\eqref{eqn_secAlgebraSM_PregularSplitCond_OurSuffConds} we would require
\begin{equation}\label{eqn_secAlgebraSM_PregularSplitCond_CoarseSuffConds}
\lambda_{min}( \mathcal{A}_i ) \geq 2 \rho(\mathtt{F}_i).
\end{equation}
\noindent Since we still have the entry-wise comparison of $u_{\ell} |\mathtt{F}_i|$ and $|\mathcal{A}_i|$, see~\eqref{eqn_secAlgebraSM_Mmtrx_mathttFi_leq_usmathcalAi}, and we have knowledge of the sign distribution of the entries of these matrices, we could arrive at some comparison theorem for $\rho(\mathtt{F}_i)$ and $\rho(\mathcal{A}_i)$, so that~\eqref{eqn_secAlgebraSM_PregularSplitCond_CoarseSuffConds} would relate the condition number $\rho(\mathcal{A}_i) / \lambda_{min}( \mathcal{A}_i )$ with the used precision $u_{\ell}$, obtaining the type of condition we encounter in~\cite[Section 3.2.2, equation (14)]{schneck2021impact} or~\cite[Section 5]{anzt2019adaptive}. The above derivation illustrates that our results are more nuanced compared to the existing ones, also in treating Schwarz methods (and their convergence) as standalone methods. 
We focus on more particular systems (in the sense of the $M$-matrix property, which together with symmetry constitutes a subclass of SPD matrices) compared to the existing literature and we carefully exploit this extra information by the specialized rounding techniques.

Next, we show results analogous to the experiments considered in Section~\ref{sec_AlgebraSM_NonSymmCase}. 
We consider the same problem as in~\eqref{eqn_secAlgebraSM_Mmtrx_Lu_eq_f}--\eqref{eqn_secAlgebraSM_Mmtrx_Lu_eq_f_DefinitionOf_L} but omit the advection terms so that we obtain Stieltjes matrices after discretization, using the following parameters.
\paragraph*{Problem 4} (analogue of Problem 1)
\begin{equation*}
\eta(\mathbf{x}) := x_1^2 \cos(x_1+x_2)^2, \;\, \alpha(\mathbf{x}) := (x_1+x_2)^2 e^{x_1-x_2}
\quad \mathrm{and} \quad
b_1(\mathbf{x}) =  b_2(\mathbf{x}) = 0.
\end{equation*}
\paragraph*{Problem 5} (analogue of Problem 2)
\begin{equation*}
\eta(\mathbf{x}) := 500 x_1 + x_2, \;\, \alpha(\mathbf{x}) := 1 + 9(x_1+x_2)
\quad \mathrm{and} \quad
b_1(\mathbf{x}) =  b_2(\mathbf{x}) = 0.
\end{equation*}
\paragraph*{Problem 6} (analogue of Problem 3)
\begin{equation*}
\eta \equiv 0, \;\, 
\alpha (\mathbf{x})= \left\{
\begin{aligned}
& 10^6 \quad \mathrm{if} \; \|\mathbf{x}-[0.5 \; 0.1]^T\| < 0.25, \\
& 1 \quad \mathrm{otherwise},
\end{aligned} \right.
\quad \mathrm{and} \quad
b_1(\mathbf{x}) = b_2(\mathbf{x}) = 0. 
\end{equation*}
\noindent We note that for Problem 4 we added non-constant reaction and diffusion coefficients (otherwise omitting the advection term leads to the standard Poisson problem).

\begin{figure}[h!]
\centering
\resizebox{1\textwidth}{!}{
	\includegraphics{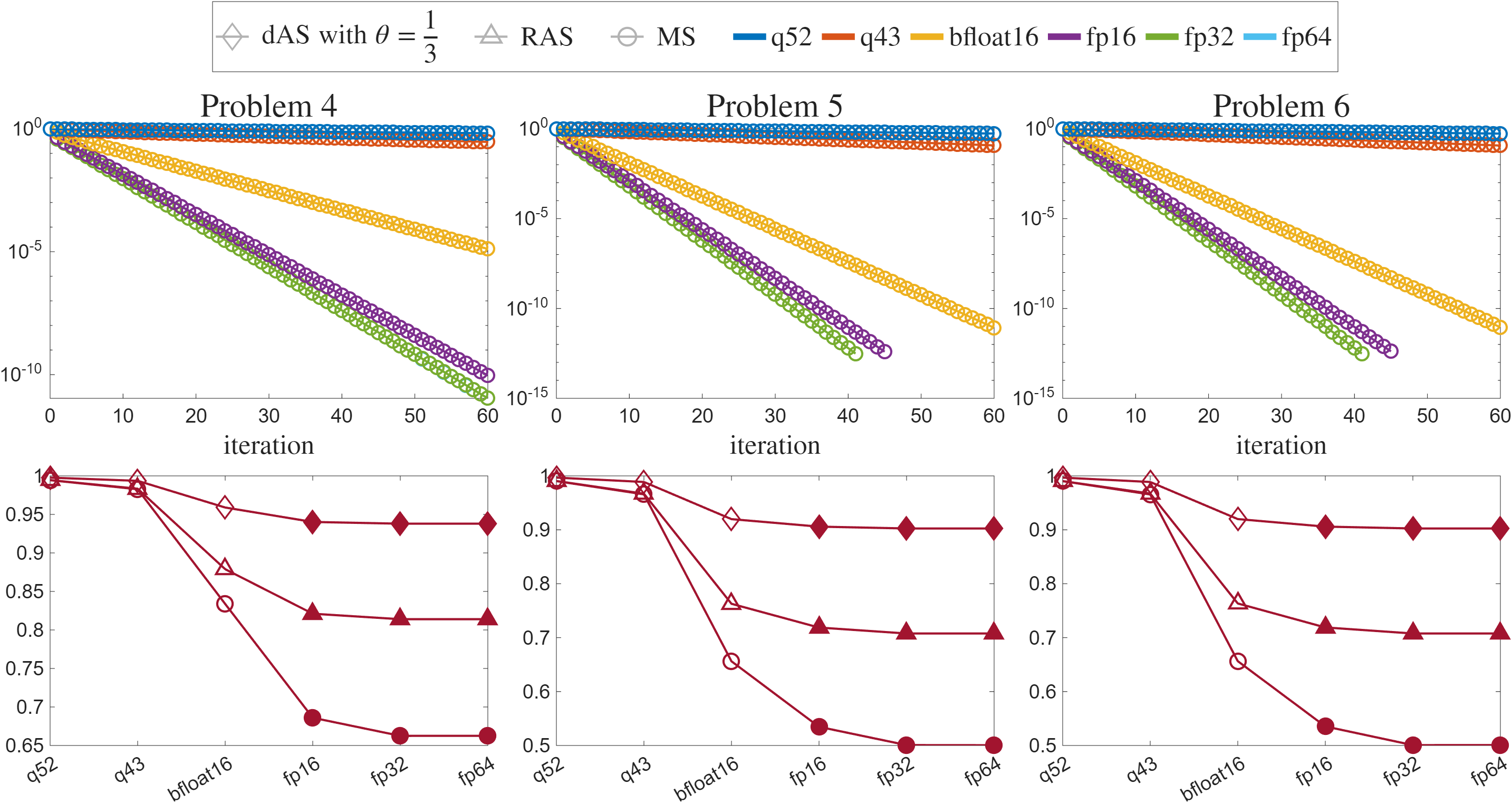}
}
\caption{Top: the 2-norm of the error of the multiplicative Schwarz method for different choices of $u_{\ell}$ and $u_w$ corresponding to \texttt{fp64}, using the \texttt{chop} toolbox. For \texttt{fp32} and \texttt{fp64} the graphs are indiscernible from each other. Bottom: the observed convergence factor $\rho_{\mathrm{conv}}$ for all Schwarz methods and choices of $u_{\ell}$; if the conditions~\eqref{eqn_secAlgebraSM_mathcalAiInvFi_leq_1},~\eqref{eqn_secAlgebraSM_Mmtrx_mathcalAiInv_geq_mathcalAiInvFimathcalAiInv} and~\eqref{eqn_secAlgebraSM_PregularSplitCond_OurSuffConds} are satisfied for both subdomain $i=1,2$ for a certain $u_{\ell}$, then the marker is filled. For example, for Problem 4 the conditions~\eqref{eqn_secAlgebraSM_mathcalAiInvFi_leq_1},~\eqref{eqn_secAlgebraSM_Mmtrx_mathcalAiInv_geq_mathcalAiInvFimathcalAiInv}  and~\eqref{eqn_secAlgebraSM_PregularSplitCond_OurSuffConds} are satisfied for both $i=1,2$ starting from \texttt{fp16}.}\label{fig_secAlgebraSM_Stieltjes_ConvCrvs_chop}
\end{figure}

\begin{figure}[h!]
\centering
\resizebox{1\textwidth}{!}{
	\includegraphics{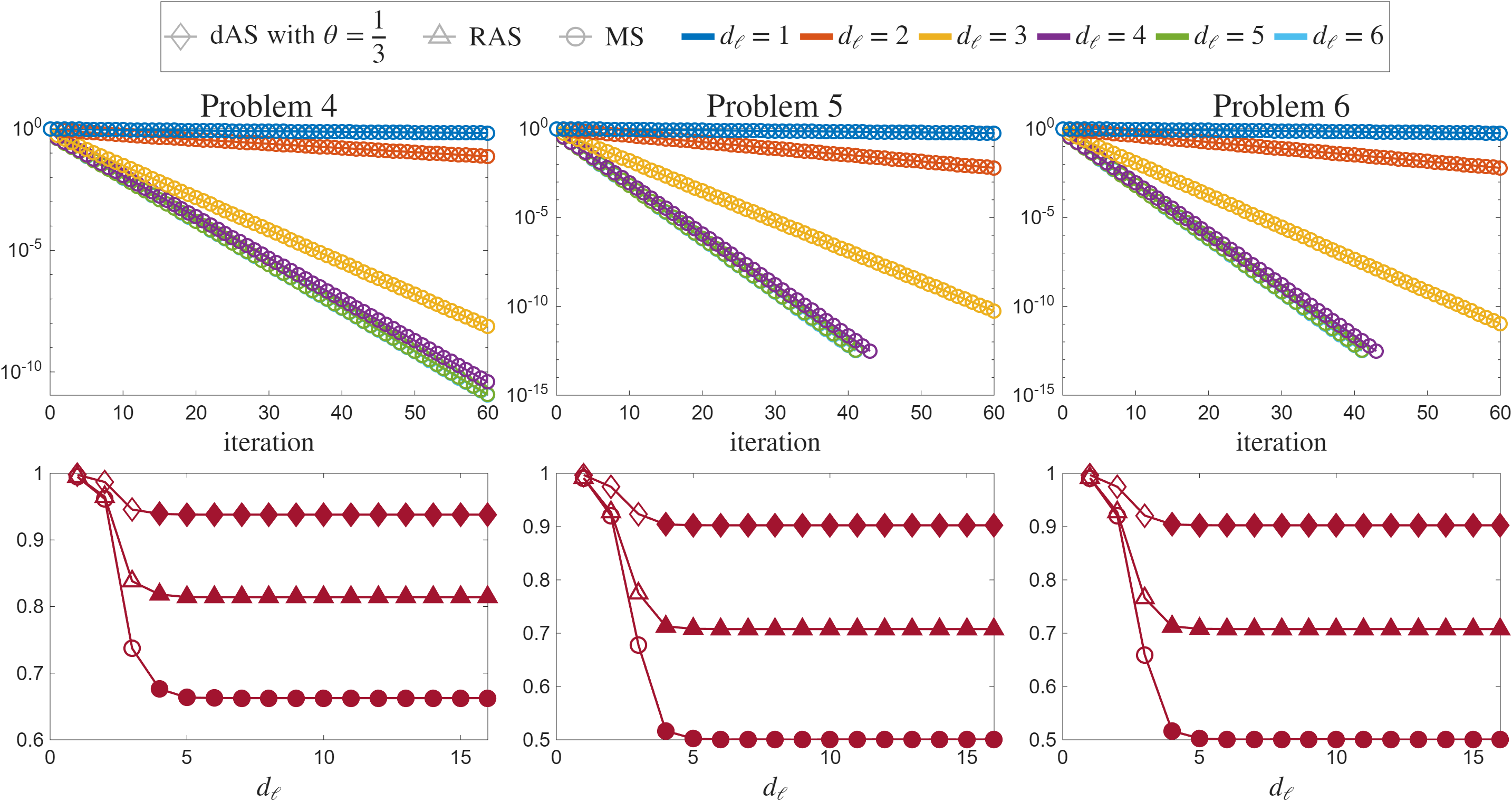}
}
\caption{Top: the 2-norm of the error of the multiplicative Schwarz method for different choices of $d_{\ell}$ and $u_w$ corresponding to \texttt{fp64}, using the \texttt{advanpix} toolbox. Bottom: the observed convergence factor $\rho_{\mathrm{conv}}$ for different methods and choices of $d_{\ell}$; if the conditions~\eqref{eqn_secAlgebraSM_mathcalAiInvFi_leq_1},~\eqref{eqn_secAlgebraSM_Mmtrx_mathcalAiInv_geq_mathcalAiInvFimathcalAiInv} and~\eqref{eqn_secAlgebraSM_PregularSplitCond_OurSuffConds} are satisfied for both subdomains $i=1,2$ for a certain $d_{\ell}$, then the marker is filled.}\label{fig_secAlgebraSM_Stieltjes_ConvCrvs_advanpix}
\end{figure}

The same questions as before are of interest\footnote{And the same caveat of the loss of efficiency with routine and naive swapping of rounding mode, see page 10, below Theorem~\ref{prop_secAlgebraSM_ASRASMSconv_MmtrxCase}.}. We fix $N=50^2$ and show the convergence curves and the observed convergence factor $\rho_{\mathrm{conv}}$ in Figure~\ref{fig_secAlgebraSM_Stieltjes_ConvCrvs_chop} (using the \texttt{chop} toolbox) and Figure~\ref{fig_secAlgebraSM_Stieltjes_ConvCrvs_advanpix} (using \texttt{advanpix} toolbox). We draw very similar conclusions to the ones in Section~\ref{sec_AlgebraSM_NonSymmCase}. We note that the additional condition~\eqref{eqn_secAlgebraSM_PregularSplitCond_OurSuffConds} was almost always weaker than~\eqref{eqn_secAlgebraSM_Mmtrx_mathcalAiInv_geq_mathcalAiInvFimathcalAiInv} and comparable to~\eqref{eqn_secAlgebraSM_mathcalAiInvFi_leq_1}. However, just as in the non-symmetric case in Section~\ref{sec_AlgebraSM_NonSymmCase}, the differences were small (e.g., for \texttt{fp16} and \texttt{bfloat16} for \texttt{chop} or for neighboring precisions for \texttt{advanpix}).

\begin{figure}[h!]
\centering
\resizebox{.9\textwidth}{!}{
	\includegraphics{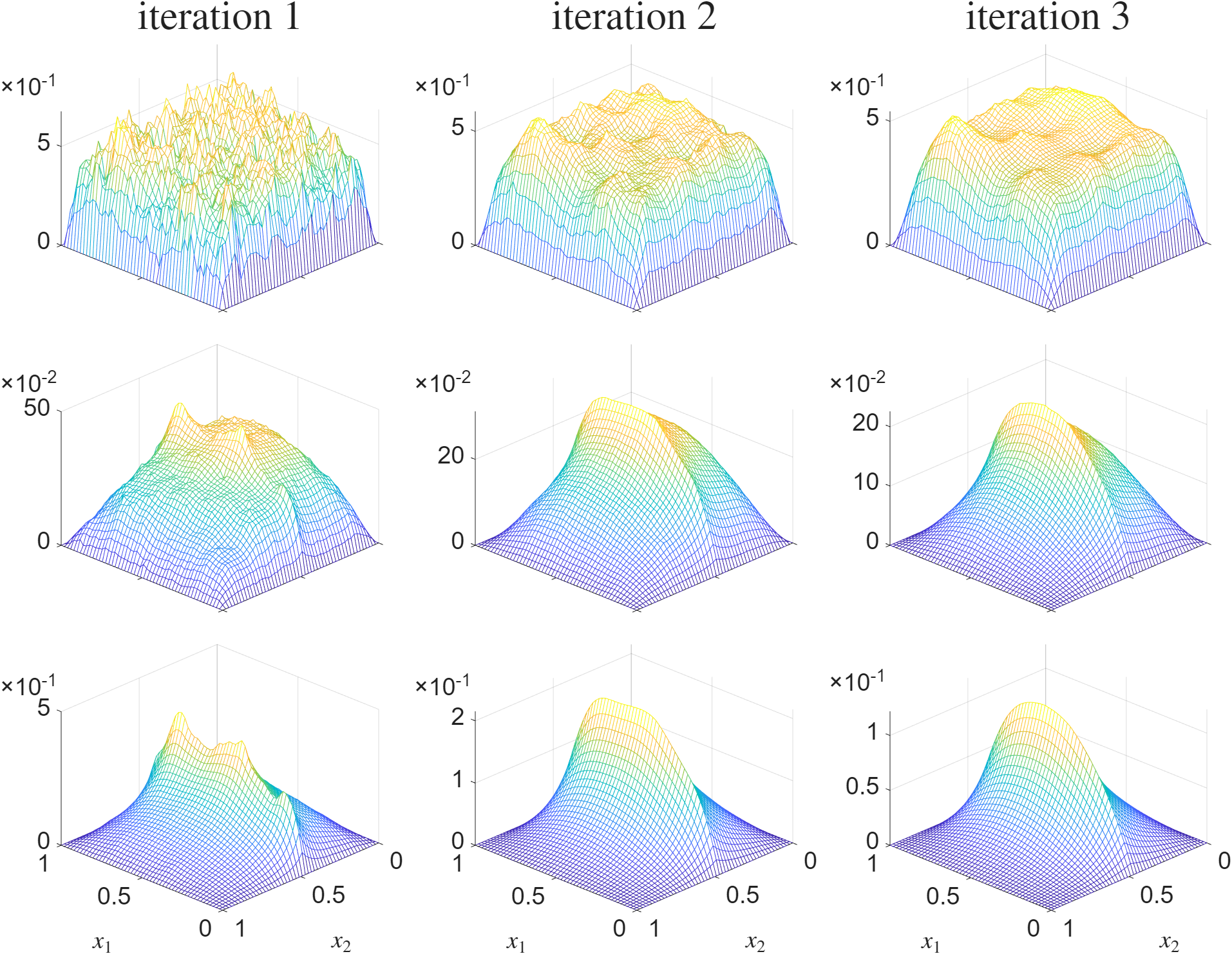}
}
\caption{The errors of the multiplicative Schwarz method used for Problem 5 with $N=2500$ during the initial iterations using \texttt{q43} (top), \texttt{bfloat16} and (middle) \texttt{fp16} (bottom). Up to scaling, the error for \texttt{q52} is analogous to the top row and the errors for \texttt{fp32} and \texttt{fp64} are are analogous to the bottom row.}\label{fig_secAlgebraSM_Stieltjes_ErrPlts}
\end{figure}

In Figure~\ref{fig_secAlgebraSM_Stieltjes_ErrPlts} we plot the error of the multiplicative Schwarz method at iterations 1, 2 and~3 for different choices of the precision $u_{\ell}$ for Problem 5 and see, generally speaking, similar results to Figure~\ref{fig_secAlgebraSM_Mmtrx_ErrPlts_chop}. Our experience with dAS and RAS is fairly similar.

\begin{figure}[t]
\centering
\resizebox{.9\textwidth}{!}{
	\includegraphics{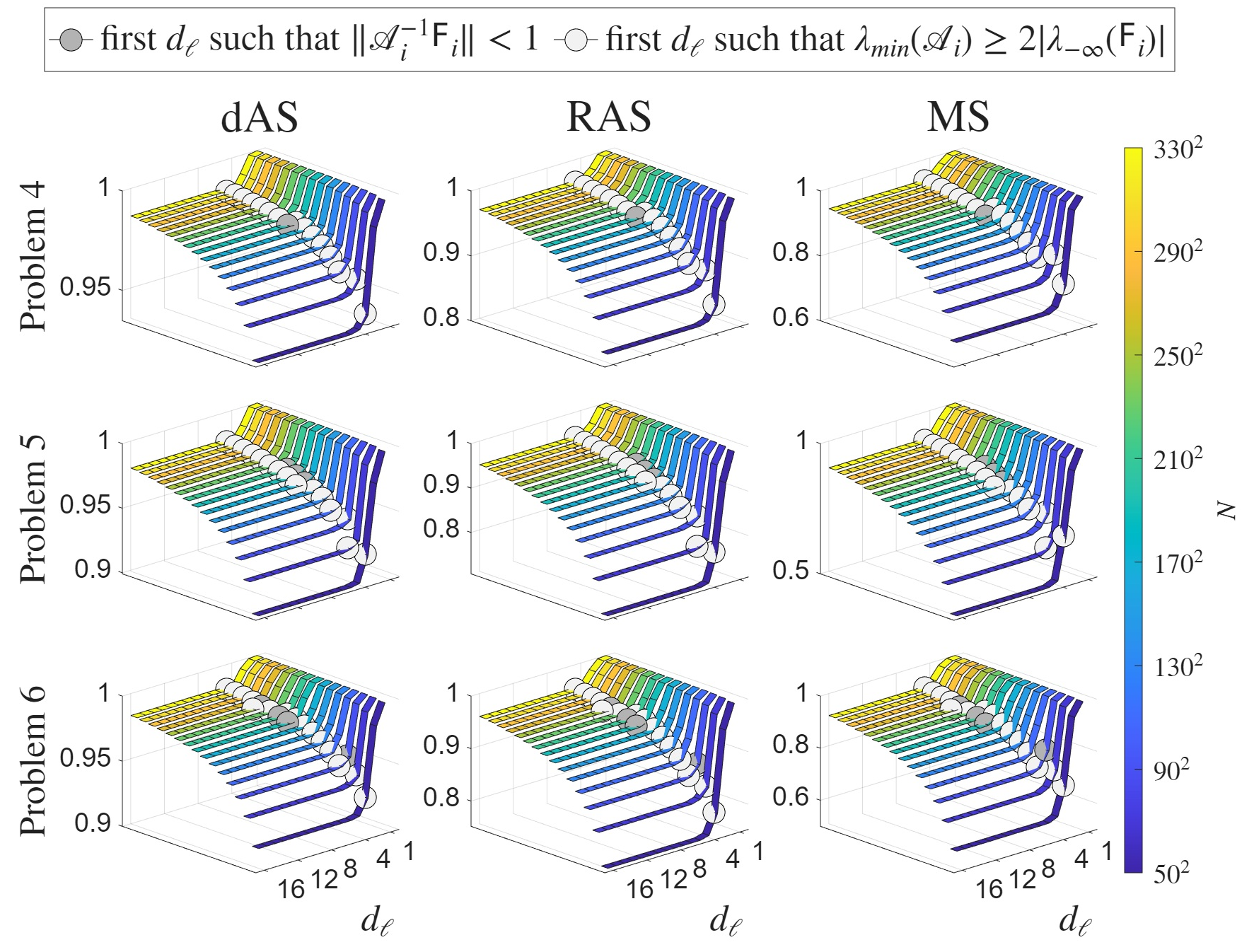}
}
\caption{We show $\rho_{\mathrm{conv}}$ of the Schwarz method for $d_{\ell}=1,2,\dotsc ,16$ and different problem sizes $N$. We highlight the first $d_{\ell}$ for which the conditions~\eqref{eqn_secAlgebraSM_mathcalAiInvFi_leq_1} and~\eqref{eqn_secAlgebraSM_PregularSplitCond_OurSuffConds} 
are satisfied by \tikzcircle[black,fill=lightgray]{3pt} and \tikzcircle[black]{2.8pt}, respectively.}\label{fig_secAlgebraSM_Stieltjes_ConvFacts}
\end{figure}

Looking at the observed convergence factors $\rho_{\mathrm{conv}}$ in Figure~\ref{fig_secAlgebraSM_Stieltjes_ConvFacts}, similarly to the non-symmetric case, the dominant eigenmodes of $T_{\star}$ with $\star \in \{ \mathrm{AS}, \mathrm{RAS}, \mathrm{MS} \}$ appear to be well-captured already with limited precision, e.g., $d_{\ell}=4\sim 6$, \emph{and} the other eigenmodes are not too sensitive with respect to small perturbations of the subdomain solves and stay ``non-dominant''. 

The theoretical results only hold if all of the conditions~\eqref{eqn_secAlgebraSM_mathcalAiInvFi_leq_1},~\eqref{eqn_secAlgebraSM_Mmtrx_mathcalAiInv_geq_mathcalAiInvFimathcalAiInv} and~\eqref{eqn_secAlgebraSM_PregularSplitCond_OurSuffConds} hold true but the model problems suggest that either of the conditions~\eqref{eqn_secAlgebraSM_mathcalAiInvFi_leq_1} or~\eqref{eqn_secAlgebraSM_PregularSplitCond_OurSuffConds} give a good indicator. However, we note that the condition~\eqref{eqn_secAlgebraSM_PregularSplitCond_OurSuffConds} becomes \emph{much} more pessimistic if we omit the re-scaling. For example, if we take $D_i^{(r,c)}=I_{N_i}$ for Problem 6, then~\eqref{eqn_secAlgebraSM_PregularSplitCond_OurSuffConds} is satisfied only for $d_{\ell} \gtrapprox 10$, i.e., long after the convergence factor has in fact stabilized at the final value. The same is true if we replace the condition~\eqref{eqn_secAlgebraSM_PregularSplitCond_OurSuffConds} with a cruder version relating to the condition number of $\mathcal{A}_i$, see~\eqref{eqn_secAlgebraSM_PregularSplitCond_CoarseSuffConds} and below. The condition~\eqref{eqn_secAlgebraSM_mathcalAiInvFi_leq_1}, however, has been fairly robust, localizing fairly accurately the optimal $d_{\ell}$ regardless of the employed scaling. Also, similarly to the non-symmetric case, the convergence factor graph becomes notably ``less smooth'' but otherwise qualitatively similar. The ``non-smoothness'' of the convergence factor for Problem 6 and the smallest mesh resolution, i.e., $N=50^2$ also stands out. The reason is not due to the low-precision use -- the algorithm has simply essentially converged after 60 iterations as we have $0.55^{60}\approx 2.6\times 10^{-16}$; this is also easy to check by inspecting the error plots directly.

\begin{figure}[t]
\centering
\resizebox{.9\textwidth}{!}{
	\includegraphics{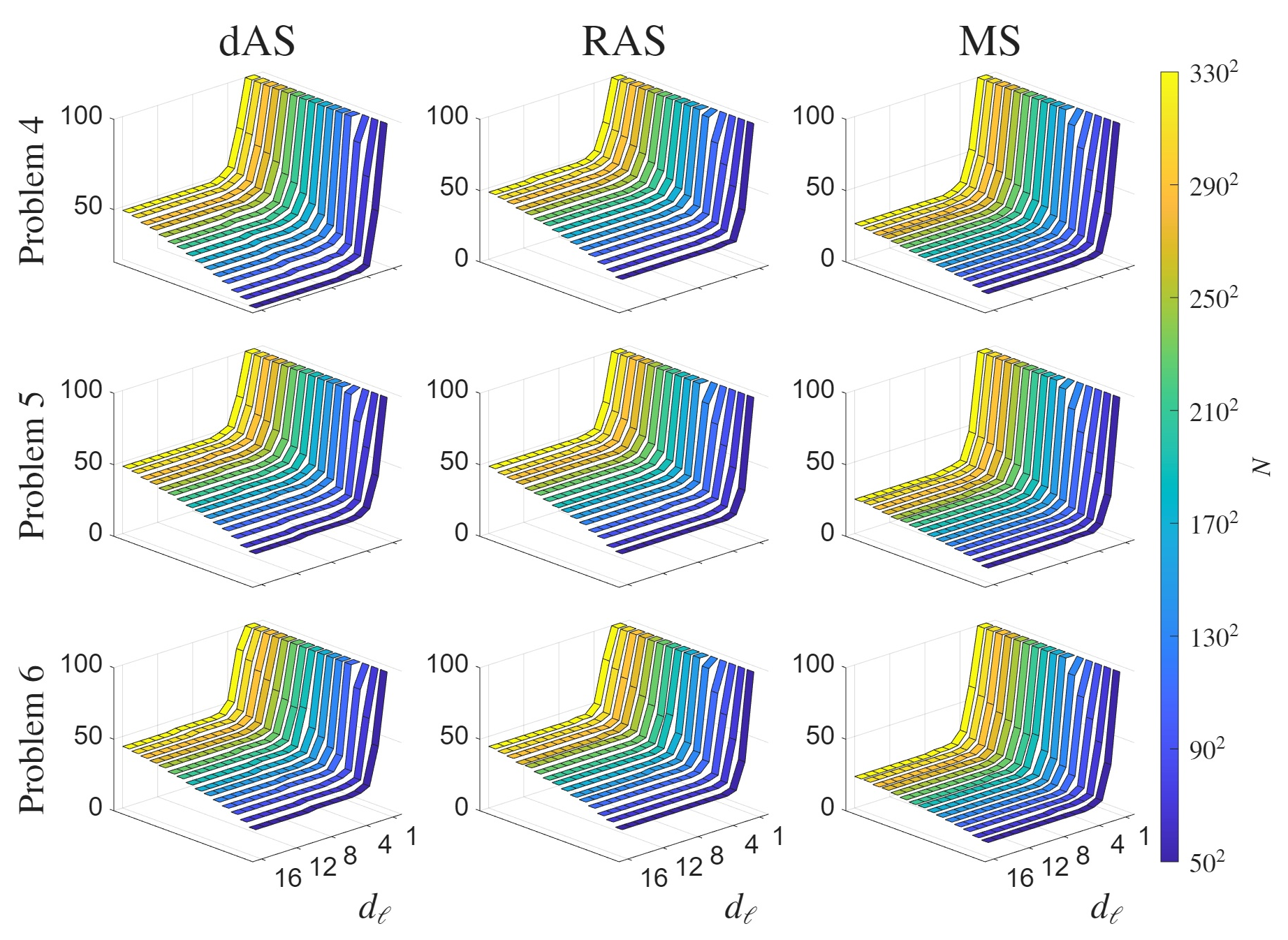}
}
\caption{The number of preconditioned GMRES iterations to reduce the relative residual below $10^{-12}$ capped at 100.}\label{fig_secAlgebraSM_Stieltjes_NmbGMRESIter}
\end{figure}

Last, we use preconditioned GMRES with multiprecision dAS\footnote{In practice, we would take advantage of the symmetry of the dAS as a preconditioner and would run a left-preconditioned CG. Here we use GMRES simply to keep the preconditioner results comparable to all of the other methods.}, RAS and MS as the left preconditioners and with relative residual tolerance $10^{-12}$, zero initial 
approximation and maximum number of iterations set to 100. We show the GMRES convergence curves and the number of iterations in 
Figure~\ref{fig_secAlgebraSM_Stieltjes_NmbGMRESIter}. We see that the number of iterations again stays mostly stable with respect to changing $d_{\ell}$ for a fixed $N$ and, moreover, the conditions~\eqref{eqn_secAlgebraSM_mathcalAiInvFi_leq_1} or~\eqref{eqn_secAlgebraSM_PregularSplitCond_OurSuffConds} still work as a reasonably accurate indicator for the choice of the number of digits $d_{\ell}$.


\section{Comparison of the standard and the multiprecision additive Schwarz methods}\label{sec_Comp_fpSM_mpSM}

So far we have studied multiprecision Schwarz methods as solvers and provided convergence conditions based on the iteration operator of the multiprecision Schwarz methods. In this section we brief\/ly explore the natural complementary approach, i.e., interpretation of the multiprecision Schwarz methods as modifications of the full-precision methods, using perturbation theory. We present the idea only for the additive Schwarz methods and leave the detailed analysis (also including other methods) as an open problem for future research.

We start by assuming the set-up of Section~\ref{sec_AlgebraSM_NonSymmCase}, namely, similarly to~\eqref{eqn_secAlgebraSM_tildemathcalAiInvFi_FormulaForNeumannExpansion} we write
\begin{equation*}
\tilde{A}_i^{-1} - A_i^{-1} = 
\mu D^{(c)}_i \left( \tilde{\mathcal{A}}^{-1}_i - \mathcal{A}_i^{-1} \right) D^{(r)}_i = 
\mu D^{(c)}_i \left( (I+\mathcal{A}_i^{-1}\mathtt{F}_i)^{-1} - I \right) \mathcal{A}_i^{-1} D^{(r)}_i ,
\end{equation*}
\noindent and denoting $\mathcal{E}_i := D^{(c)}_i \left( (I+\mathcal{A}_i^{-1}\mathtt{F}_i)^{-1} - I \right) (D^{(c)}_i)^{-1}$ we obtain
\begin{equation}
\tilde{A}_i^{-1} - A_i^{-1} = \mathcal{E}_i A_i^{-1}.
\end{equation}
\noindent As a result, if we assume\footnote{This assumption is analogous to~\eqref{eqn_secAlgebraSM_mathcalAiInvFi_leq_1}. In fact, the derivations requiring~\eqref{eqn_secAlgebraSM_mathcalAiInvFi_leq_1} can be carried out analogously even if we assume~\eqref{eqn_secCnvrgncBnds_invAimathttFi_leq_eps} instead of~\eqref{eqn_secAlgebraSM_mathcalAiInvFi_leq_1} but the derivation becomes somewhat more lengthy.}
\begin{equation}\label{eqn_secCnvrgncBnds_invAimathttFi_leq_eps}
\| \mathcal{A}_i^{-1}\mathtt{F}_i \| \leq \epsilon < \frac{1}{2},
\end{equation}
\noindent for some $\epsilon \in (0,1/2)$, then 
\begin{equation*}
\| \mathcal{E}_i \|  = \left\|  D^{(c)}_i \left( \sum\limits_{k=1}^{\infty} (-1)^{k} \left( \mathcal{A}_i^{-1}\mathtt{F}_i \right)^{k} \right) (D^{(c)}_i)^{-1} \right\| \leq \kappa \left( D_i^{(c)} \right) \epsilon \frac{1}{1-\epsilon} < 2\epsilon \kappa \left( D_i^{(c)} \right),
\end{equation*}
\noindent and thus
\begin{equation}
\| \tilde{A}_i^{-1} - A_i^{-1} \| \leq 2\epsilon \kappa \left( D_i^{(c)} \right) \| A_i^{-1} \|,
\end{equation}
\noindent where $\kappa(\cdot)$ denotes the condition number with respect to the norm $\|\cdot\|$. In words, a (small) perturbation of the scaled subdomain matrices can perturb the subdomain solves proportionally to the given subdomain solve norm and to the condition number of the column-scaling matrix. Hence we hope to get a well-scaled subdomain matrices, for which the scaling can be mostly done by row-scaling (recall that, conveniently, the row-scaling takes precedence in~\cite[Algorithms 2.3 and 2.4]{higham2019squeezing}). In the case of (damped) AS, one of the main advantages of (damped) AS is the symmetry and hence symmetrical scaling is of special interest, where we would also have to hope for a well-scaled matrix. Also, notice that~\eqref{eqn_secCnvrgncBnds_invAimathttFi_leq_eps} is similar to the assumption~\cite[equation (12)]{schneck2021impact} but for the inverses and after the scaling. 

Next, we insert~\eqref{eqn_secCnvrgncBnds_invAimathttFi_leq_eps} into the definition of the additive Schwarz methods in~\eqref{eqn_secIntro_MultSchwrz_SchwarzMthds_DefOfIterMtrcsAndPreconds} and get
\begin{equation}\label{eqn_secCnvrgncBnds_MultPrecisionPrecSysAsPerturbOfFullPrecisionPrecSys}
\tilde{M}_{\star}^{-1}A = M_{\star}^{-1}A +  E_{\star}
\quad \mathrm{where} \quad
E_{\star}  := \begin{cases}
\theta \sum\limits_{i=1}^{p} R_i^T \mathcal{E}_i A_i^{-1} R_iA, &\mathrm{for \; (damped) \;AS}, \\
\sum\limits_{i=1}^{p} \overline{R}_i^T \mathcal{E}_i A_i^{-1} R_iA, &\mathrm{for \; RAS}.
\end{cases}
\end{equation}
\noindent Recalling the matrix definitions in~\eqref{eqn_secIntro_SM__Ri_eq_IdZeroPi} and below, we can write
\begin{equation*}
A_i^{-1} R_iA = A_i^{-1} \begin{bmatrix} I_{N_i} & 0 \end{bmatrix} \Pi_{i} \Pi_{i}^T
\begin{bmatrix} A_i & K_{i} \\ L_{i} & A_{\neg i} \end{bmatrix} \Pi_{i} = 
\begin{bmatrix} I_{N_i} & A_i^{-1}K_{i} \end{bmatrix} \Pi_{i}.
\end{equation*}
\noindent and hence 
\begin{equation*}
E_{\star} = 
\begin{cases}
AS,\theta \, : \quad  &\theta \sum\limits_{i=1}^{p} \Pi_i^T \begin{bmatrix} \mathcal{E}_i & \mathcal{E}_i A_i^{-1}K_{i} \\ 0 & 0 \end{bmatrix} \Pi_{i}, \\
RAS \, : \quad  &\sum\limits_{i=1}^{p} \overline{\Pi}_i^T \begin{bmatrix} \left(\mathcal{E}_i\right)_{1:\overline{N}_i,:} & \left( \mathcal{E}_i A_i^{-1}K_{i} \right)_{1:\overline{N}_i,:} \\ 0 & 0 \end{bmatrix} \Pi_{i}.
\end{cases}
\end{equation*}
\noindent In fact, we can further rewrite this as
\begin{equation}\label{eqn_secCnvrgncBnds_ErrMtrx_E}
E_{\star} = 
\begin{cases}
AS,\theta \, : \quad  &\theta \sum\limits_{i=1}^{p} \Pi_i^T \begin{bmatrix} \mathcal{E}_i & 0 \\ 0 & 0 \end{bmatrix} \begin{bmatrix} I_{N_i} & A_i^{-1}K_{i} \\ 0 & 0 \end{bmatrix} \Pi_{i}, \\
RAS \, : \quad  &\sum\limits_{i=1}^{p} \overline{\Pi}_i^T 
\begin{bmatrix} \begin{bmatrix}
( \mathcal{E}_i )_{1:\overline{N}_i,:} \\ 0_{\overline{N}_i+1:N_i,:}
\end{bmatrix} & 0 \\ 0 & 0 \end{bmatrix}
\begin{bmatrix} I_{N_i} & A_i^{-1}K_{i}  \\ 0 & 0 \end{bmatrix} \Pi_{i},
\end{cases}
\end{equation}
\noindent where we clearly see the main difference between the two methods in the ``double-counting of the overlap'', see~\cite{efstathiou2003restricted,frommer2001algebraic}. We outline the importance of such analysis in the following remarks.

\begin{remark}
In our opinion, a more detailed study of the \emph{interaction} of the error introduced by low-precision subdomain solves and the structure of the Schwarz methods can lead to more insightful analysis. In this way, we would be able to obtain better guidance for choosing the precision $u_{\ell}$.
\end{remark}

\begin{remark}
The formulation~\eqref{eqn_secCnvrgncBnds_MultPrecisionPrecSysAsPerturbOfFullPrecisionPrecSys} allows us to use results from the convergence theory of Krylov subspace methods for perturbed systems. In particular, results like~\cite[Theorem~1.1 and Corollary~1.2]{blechta2021stability} or~\cite[Theorems 2.1 and 2.3]{embree2013gmres} give direct recipes for bounding the convergence delay of GMRES due to the subdomain solves in the precision $u_{\ell}$, assuming $E_{\star}$ is in some sense small. Perturbation results exist also for CG and can be of use for the (damped) AS, see the seminal papers~\cite{greenbaum1989behavior,strakos1991real} and also~\cite[Section 5.9]{liesen2013krylov} for further references.
\end{remark}

\section{Conclusion and future work}\label{sec_Conclusion}
We have proposed and analyzed multiprecision Schwarz methods that are specifically tailored for problems where we can guarantee the methods' convergence -- problems where the system matrix is an $M$-matrix. Using specific rounding techniques, we were able to preserve the convergence property and suggest several natural conditions for choosing a suitable precision depending on the problem. We presented several numerical experiments on PDE model problems that support our theoretical results and further illustrate that our proposed conditions describe the ideal precision for the subdomain solves. As future work we intend to consider generalizations for multiple subdomains and/or ``interface conditions'' in the sense of~\cite{gander2012optimal}.

An understanding of the interaction of the subdomain matrices $A_i^{-1}K_i$ and the matrices $\mathcal{E}_i$ for all three classical Schwarz methods for a wider variety of problems would be certainly interesting and we leave it open as a possibility for future research. Also, it has been shown that it is often more suitable to bound the norm of the iteration matrix (and hence of the error) over two or more iterations due to the nature of the underlying PDE analysis, see, e.g.,~\cite{outrata2022onmras,gander2012optimal}. Exploiting this to get a better grasp on the multiprecision Schwarz methods as stand-alone solvers would be useful. Naturally, extending this analysis to preconditioning or rather understanding how to do that would be also of clear interest. This would be likely overlapping with the so-called double-sweeping preconditioners and their analysis.

Finally, the authors would like to thank Petr Vacek for many stimulating and fruitful discussion on the topic of multiprecision computations, Theo Mary for enlighting discussions which led to our
Remark~\ref{Theo:remark},  and Andreas Frommer for insightful feedback on the submitted manuscript, leading to its notable improvement.

\bibliographystyle{plainnat}
\bibliography{biblio}

\end{document}